\documentstyle[amssymb,amsmath]{article}
\newtheorem{th}{Theorem}[section]
\newtheorem{lem}[th]{Lemma}

\newtheorem{cor}[th]{Corollary}
\newtheorem{defn}[th]{Definition}
\newenvironment{defn-new}{\begin{defn} \em}{\end{defn}}
\newtheorem{rem}[th]{Remark}
\newenvironment{rem-new}{\begin{rem} \em}{\end{rem}}
\newtheorem{ex}[th]{Example}
\newenvironment{ex-new}{\begin{ex} \em}{\end{ex}}

\newtheorem{notation}[th]{Notation}
\newenvironment{notation-new}{\begin{rem} \em}{\end{rem}}

\newenvironment{agr-new}{\begin{rem} \em}{\end{rem}}

\makeatletter \@addtoreset{equation}{section} \makeatother

\makeatletter \@addtoreset{figure}{section} \makeatother

\begin{document}

\begin{center}
{\LARGE {\bf On $(N(k),\xi )$-semi-Riemannian manifolds: Semisymmetries}}\\[%
0pt]
\bigskip \bigskip {\large {\bf Mukut Mani Tripathi and Punam Gupta}}
\end{center}

{\bf Abstract.} $(N(k),\xi )$-semi-Riemannian manifolds are defined.
Examples and properties of $(N(k),\xi )$-semi-Riemannian manifolds are
given. Some relations involving ${\cal T}_{\!a}$-curvature tensor in $%
(N(k),\xi )$-semi-Riemannian manifolds are proved. $\xi $-${\cal T}_{\!a}$%
-flat $(N(k),\xi )$-semi-Riemannian manifolds are defined. It is proved that
if $M$ is an $n$-dimensional $\xi $-${\cal T}_{\!a}$-flat $(N(k),\xi )$%
-semi-Riemannian manifold, then it is $\eta $-Einstein under an algebraic
condition. We prove that a semi-Riemannian manifold, which is $T$-recurrent
or $T$-symmetric, is always $T$-semisymmetric, where $T$ is any tensor of
type $(1,3)$. $\left( {\cal T}_{\!a}, {\cal T}_{\!b}\right) $-semisymmetric
semi-Riemannian manifold is defined and studied. The results for ${\cal T}%
_{\!a}$-semisymmetric, ${\cal T}_{\!a}$-symmetric, ${\cal T}_{\!a}$%
-recurrent $(N(k),\xi )$-semi-Riemannian manifolds are obtained. The
definition of $({\cal T}_{\!a},S_{{\cal T}_{b}})$-semisymmetric
semi-Riemannian manifold is given. $({\cal T}_{\!a},S_{{\cal T}_{b}})$%
-semisymmetric $(N(k),\xi )$-semi-Riemannian manifolds are classified. Some
results for ${\cal T}_{\!a}$-Ricci-semisymmetric $(N(k),\xi )$%
-semi-Riemannian manifolds are obtained. \medskip

\medskip

\noindent {\bf Keywords.} $(N(k),\xi )$-semi-Riemannian manifold; ${\cal T}$%
-curvature tensor; quasi-conformal curvature tensor; conformal curvature
tensor; conharmonic curvature tensor; concircular curvature tensor;
pseudo-projective curvature tensor; projective curvature tensor; ${\cal M}$%
-projective curvature tensor; ${\cal W}_{i}$-curvature tensors $(i=0,\ldots
,9)$; ${\cal W}_{j}^{\ast }$-curvature tensors $(j=0,1)$; $\eta $-Einstein
manifold; Einstein manifold; $N(k)$-contact metric manifold; $\left(
\varepsilon \right) $-Sasakian manifold; Sasakian manifold; Kenmotsu
manifold; $\left( \varepsilon \right) $-para-Sasakian manifold;
para-Sasakian manifold; $(N(k),\xi )$-semi-Riemannian manifolds; $\left( 
{\cal T}_{\!a},{\cal T}_{\!b}\right) $-semisymmetric semi-Riemannian
manifold; $({\cal T}_{\!a},S_{{\cal T}_{b}})$-semisymmetric semi-Riemannian
manifold and $\xi $-${\cal T}_{\!a}$-flat $(N(k),\xi )$-semi-Riemannian
manifold.

\section{Introduction}

Let $(M,g)$ be an $n$-dimensional semi-Riemannian manifold and ${\frak X}(M)$
the Lie algebra of vector fields in $M$. Throughout the paper we assume that 
$X,Y,Z,U,V,W\in {\frak X}(M)$, unless specifically stated otherwise.\medskip

A semi-Riemannian manifold $M$ is said to be flat if $R(X,Y)Z=0$. It is said
to be $\xi $-flat if $R(X,Y)\xi =0$, where $\xi $ is a non-null unit vector
field in $M$. The condition of $\xi $-flatness is weaker than the condition
of flatness. In 2006, De and Biswas \cite{De-Biswas-06} studied the $\xi $%
-conformally flat contact metric manifolds with $\xi \in N(k)$. They proved
that a contact metric manifold with $\xi \in N(k)$ is $\xi $-conformally
flat if and only if it is $\eta $-Einstein manifold. Recently, in 2010,
Dwivedi and Kim \cite{Dwivedi-Kim-10} proved that a Sasakian manifold is $%
\xi $-conharmonically flat if and only if it is $\eta $-Einstein.\medskip

A semi-Riemannian manifold $M$ is said to be semisymmetric \cite{Szabo-82}\
if it satisfies $R(X,Y)\cdot R=0$, where $R\left( X,Y\right) $ acts as a
derivation on $R$. Semisymmetric manifold is a generalization of manifold of
constant curvature and symmetric manifold ($\nabla R=0$). A semi-Riemannian
manifold is said to be recurrent \cite{Walker-50} if it satisfies $\nabla
R=\alpha \otimes R$, where $\alpha $ is $1$-form. In 1972, Takagi \cite%
{Takagi-72} gave an example of Riemannian manifolds satisfying $R(X,Y)\cdot
R=0$ but not $\nabla R=0$. \medskip

A semi-Riemannian manifold $M$ is said to be Ricci-semisymmetric \cite%
{Deszcz-89} if its Ricci tensor $S$ satisfies $R(X,Y)\cdot S=0$, where $%
R\left( X,Y\right) $ acts as a derivation on $S$. Ricci-semisymmetric
manifold is a generalization of manifold of constant curvature, Einstein
manifold, Ricci symmetric manifold, symmetric manifold and semisymmetric
manifold. \medskip 

Ricci-semisymmetric manifolds are studied by Adati and Miyazawa \cite%
{Adati-Miyazawa-79}, Hong et al. \cite{Hong-Ozgur-Tripathi-06}, Pandey and
Verma \cite{Pandey-Verma-99}, Perrone \cite{Perrone-92} and Tripathi et al. 
\cite{TKYK-09}. After this \"{O}zg\"{u}r \cite{Ozgur-2005} studied the Weyl
Ricci-semisymmetric manifold. Hong, \"{O}zg\"{u}r and Tripathi (\cite%
{Hong-Ozgur-Tripathi-06}, \cite{Ozgur-Tripathi-2007}) studied the
concircular Ricci-semisymmetric manifold. \medskip

The paper is organized as follows. In Section \ref{sect-GCT}, we give the
definition of ${\cal T}$-curvature tensor. In Section~\ref{sect-NK}, we
define $(N(k),\xi )$-semi-Riemannian manifolds. Examples and properties of $%
(N(k),\xi )$-semi-Riemannian manifolds are given. $N(k)$-contact metric
manifold, $\left( \varepsilon \right) $-Sasakian, Sasakian, Kenmotsu, $%
\left( \varepsilon \right) $-para-Sasakian and para-Sasakian manifolds are
examples of $(N(k),\xi )$-semi-Riemannian manifolds. We obtain the relations
for ${\cal T}_{\!a}$-curvature tensor in $(N(k),\xi )$-semi-Riemannian
manifold. In Section~\ref{sect-T0}, the definition of $\xi $-${\cal T}_{\!a}$%
-flat $(N(k),\xi )$-semi-Riemannian manifold is given. Necessary conditions
for a $\xi $-${\cal T}_{\!a}$-flat $(N(k),\xi )$-semi-Riemannian manifold
are mentioned. It is proved that an $n$-dimensional $\xi $-${\cal T}_{\!a}$%
-flat $(N(k),\xi )$-semi-Riemannian manifold is $\eta $-Einstein under an
algebraic condition. The necessary and sufficient condition for an $n$%
-dimensional $(N(k),\xi )$-semi-Riemannian manifold to be $\xi $-${\cal T}%
_{\!a}$-flat is obtained, where ${\cal T}_{\!a}$-curvature tensor is one of
the quasi-conformal curvature tensor, conformal curvature tensor,
conharmonic curvature tensor, ${\cal M}$-projective curvature tensor or $%
{\cal W}_{2}$-curvature tensor. In Section~\ref{sect-TR}, the definition of $%
T$-recurrent, $T$-symmetric and $T$-semisymmetric semi-Riemannian manifolds
are given, where $T$ is any tensor of type $(1,3)$. It is proved that if a
semi-Riemannian manifold is $T$-recurrent or $T$-symmetric, then it is
always $T$-semisymmetric. In Section~\ref{sect-TT}, $\left( {\cal T}_{\!a},%
{\cal T}_{b}\right) $-semisymmetric semi-Riemannian manifolds are defined
and classified. It is proved that ${\cal T}_{\!a}$-semisymmetric $(N(k),\xi )
$-semi-Riemannian manifolds are either $\eta $-Einstein or Einstein and
manifold of constant curvature. ${\cal T}_{\!a}$-semisymmetric $(N(k),\xi )$%
-semi-Riemannian manifold is proved to be ${\cal T}_{\!a}$-flat under
certain condition. A ${\cal T}_{\!a}$-semisymmetric $(N(k),\xi )$%
-semi-Riemannian manifold is ${\cal T}_{\!a}$-conservative under some
condition. It is also proved that if an $(N(k),\xi )$-semi-Riemannian
manifold is of constant curvature, then it is ${\cal T}_{\!a}$-semisymmetric
under some algebraic conditions. In the last section, the definition of $(%
{\cal T}_{\!a},S_{{\cal T}_{b}})$-semisymmetric semi-Riemannian manifold is
given. $({\cal T}_{\!a},S_{{\cal T}_{b}})$-semisymmetric $(N(k),\xi )$%
-semi-Riemannian manifolds are classified. The results for ${\cal T}_{\!a}$%
-Ricci-semisymmetric $(N(k),\xi )$-semi-Riemannian manifolds are obtained.
An $n$-dimensional $(R,S_{{\cal T}_{a}})$-semisymmetric $(N(k),\xi )$%
-semi-Riemannian manifold is Einstein under an algebraic condition. An
Einstein manifold is $(R,S_{{\cal T}_{a}})$-semisymmetric. If $M$ is an
Einstein manifold such that ${\cal T}_{\!a}\in \left\{ R,{\cal C}_{\ast },%
{\cal C},{\cal L},{\cal V},{\cal M},{\cal W}_{0},{\cal W}_{0}^{\ast },{\cal W%
}_{3}\right\} $, then it is ${\cal T}_{\!a}$-Ricci-semisymmetric.

\section{${\cal T}$-curvature tensor\label{sect-GCT}}


\begin{defn-new}
\label{defn-GCT} In an $n$-dimensional semi-Riemannian manifold $\left(
M,g\right) $, {\em ${\cal T}$-curvature tensor} {\rm {\cite{Tripathi-Gupta}} 
}is a tensor of type $(1,3)$, which is defined by 
\begin{eqnarray}
{\cal T}\left( X,Y\right) Z &=&a_{0}\,R\left( X,Y\right) Z  \nonumber \\
&&+\ a_{1}\,S\left( Y,Z\right) X+a_{2}\,S\left( X,Z\right) Y+a_{3}\,S(X,Y)Z 
\nonumber \\
&&+\ a_{4}\,g\left( Y,Z\right) QX+a_{5}\,g\left( X,Z\right)
QY+a_{6}\,g(X,Y)QZ  \nonumber \\
&&+\ a_{7}\,r\left( g\left( Y,Z\right) X-g\left( X,Z\right) Y\right) ,
\label{eq-GCT}
\end{eqnarray}%
where $a_{0},\ldots ,a_{7}$ are real numbers; and $R$, $S$, $Q$ and $r$ are
the curvature tensor, the Ricci tensor, the Ricci operator and the scalar
curvature respectively.
\end{defn-new}

In particular, the ${\cal T}$-curvature tensor is reduced to

\begin{enumerate}
\item the{\em \ curvature tensor} $R$ if \vspace{-0.4cm} 
\[
a_{0}=1,\quad a_{1}=\cdots =a_{7}=0, 
\]%
\vspace{-0.5cm}

\item the {\em quasi-conformal curvature tensor} ${\cal C}_{\ast }$ \cite%
{Yano-Sawaki-68} if \vspace{-0.4cm} 
\[
a_{1}=-\,a_{2}=a_{4}=-\,a_{5},\quad a_{3}=a_{6}=0,\quad a_{7}=-\,\frac{1}{n}%
\left( \frac{a_{0}}{n-1}+2a_{1}\right) , 
\]%
\vspace{-0.5cm}

\item the {\em conformal curvature tensor} ${\cal C}$ \cite[p.~90]%
{Eisenhart-49} if \vspace{-0.4cm} 
\[
a_{0}=1,\quad a_{1}=-\,a_{2}=a_{4}=-\,a_{5}=-\,\frac{1}{n-2},\quad
a_{3}=a_{6}=0,\quad a_{7}=\frac{1}{(n-1)(n-2)}, 
\]%
\vspace{-0.5cm}

\item the {\em conharmonic curvature tensor} ${\cal L}$ \cite{Ishii-57} if 
\vspace{-0.4cm} 
\[
a_{0}=1,\quad a_{1}=-\,a_{2}=a_{4}=-\,a_{5}=-\,\frac{1}{n-2},\,\quad
a_{3}=a_{6}=0,\quad a_{7}=0, 
\]%
\vspace{-0.5cm}

\item the {\em concircular curvature tensor} ${\cal V}$ (\cite{Yano-40}, 
\cite[p. 87]{Yano-Bochner-53}) if \vspace{-0.4cm} 
\[
a_{0}=1,\quad a_{1}=a_{2}=a_{3}=a_{4}=a_{5}=a_{6}=0,\quad a_{7}=-\,\frac{1}{%
n(n-1)}, 
\]%
\vspace{-0.5cm}

\item the {\em pseudo-projective curvature tensor }${\cal P}_{\ast }$ \cite%
{Prasad-2002} if \vspace{-0.4cm} 
\[
a_{1}=-\,a_{2},\quad a_{3}=a_{4}=a_{5}=a_{6}=0,\quad a_{7}=-\,\frac{1}{n}%
\left( \frac{a_{0}}{n-1}+a_{1}\right) , 
\]%
\vspace{-0.5cm}

\item the {\em projective curvature tensor} ${\cal P}$ \cite[p. 84]%
{Yano-Bochner-53} if \vspace{-0.4cm} 
\[
a_{0}=1,\quad a_{1}=-\,a_{2}=-\,\frac{1}{(n-1)}\text{,\quad }%
a_{3}=a_{4}=a_{5}=a_{6}=a_{7}=0, 
\]%
\vspace{-0.5cm}

\item the ${\cal M}${\em -projective curvature tensor }\cite%
{Pokhariyal-Mishra-71} if \vspace{-0.4cm} 
\[
a_{0}=1,\quad a_{1}=-\,a_{2}=a_{4}=-\,a_{5}=-\frac{1}{2(n-1)},\quad
a_{3}=a_{6}=a_{7}=0, 
\]%
\vspace{-0.5cm}

\item the ${\cal W}_{0}${\em -curvature tensor} \cite[Eq. (1.4)]%
{Pokhariyal-Mishra-71} if \vspace{-0.4cm} 
\[
a_{0}=1,\quad a_{1}=-\,a_{5}=-\,\frac{1}{(n-1)},\quad
a_{2}=a_{3}=a_{4}=a_{6}=a_{7}=0, 
\]%
\vspace{-0.5cm}

\item the ${\cal W}_{0}^{\ast }${\em -curvature tensor} \cite[Eq. (2.1)]%
{Pokhariyal-Mishra-71} if \vspace{-0.4cm} 
\[
a_{0}=1,\quad a_{1}=-\,a_{5}=\frac{1}{(n-1)},\quad
a_{2}=a_{3}=a_{4}=a_{6}=a_{7}=0, 
\]%
\vspace{-0.5cm}

\item the ${\cal W}_{1}${\em -curvature tensor} \cite{Pokhariyal-Mishra-71}
if \vspace{-0.4cm} 
\[
a_{0}=1,\quad a_{1}=-\,a_{2}=\frac{1}{(n-1)},\quad
a_{3}=a_{4}=a_{5}=a_{6}=a_{7}=0, 
\]%
\vspace{-0.5cm}

\item the ${\cal W}_{1}^{\ast }${\em -curvature tensor} \cite%
{Pokhariyal-Mishra-71} if \vspace{-0.4cm} 
\[
a_{0}=1,\quad a_{1}=-\,a_{2}=-\,\frac{1}{(n-1)},\quad
a_{3}=a_{4}=a_{5}=a_{6}=a_{7}=0, 
\]%
\vspace{-0.5cm}

\item the ${\cal W}_{2}${\em -curvature tensor} \cite{Pokhariyal-Mishra-70}
if \vspace{-0.4cm} 
\[
a_{0}=1,\quad a_{4}=-\,a_{5}=-\,\frac{1}{(n-1)},\quad
a_{1}=a_{2}=a_{3}=a_{6}=a_{7}=0, 
\]%
\vspace{-0.5cm}

\item the ${\cal W}_{3}${\em -curvature tensor} \cite{Pokhariyal-Mishra-71}
if \vspace{-0.4cm} 
\[
a_{0}=1,\quad a_{2}=-\,a_{4}=-\,\frac{1}{(n-1)},\quad
a_{1}=a_{3}=a_{5}=a_{6}=a_{7}=0, 
\]%
\vspace{-0.4cm}

\item the ${\cal W}_{4}${\em -curvature tensor} \cite{Pokhariyal-Mishra-71}
if \vspace{-0.5cm} 
\[
a_{0}=1,\quad a_{5}=-\,a_{6}=\frac{1}{(n-1)},\quad
a_{1}=a_{2}=a_{3}=a_{4}=a_{7}=0, 
\]%
\vspace{-0.5cm}

\item the ${\cal W}_{5}${\em -curvature tensor} \cite{Pokhariyal-82} if 
\vspace{-0.4cm} 
\[
a_{0}=1,\quad a_{2}=-\,a_{5}=-\,\frac{1}{(n-1)},\quad
a_{1}=a_{3}=a_{4}=a_{6}=a_{7}=0, 
\]%
\vspace{-0.5cm}

\item the ${\cal W}_{6}${\em -curvature tensor} \cite{Pokhariyal-82} if 
\vspace{-0.4cm} 
\[
a_{0}=1,\quad a_{1}=-\,a_{6}=-\,\frac{1}{(n-1)},\quad
a_{2}=a_{3}=a_{4}=a_{5}=a_{7}=0, 
\]%
\vspace{-0.5cm}

\item the ${\cal W}_{7}${\em -curvature tensor} \cite{Pokhariyal-82} if 
\vspace{-0.4cm} 
\[
a_{0}=1,\quad a_{1}=-\,a_{4}=-\,\frac{1}{(n-1)},\quad
a_{2}=a_{3}=a_{5}=a_{6}=a_{7}=0, 
\]%
\vspace{-0.5cm}

\item the ${\cal W}_{8}${\em -curvature tensor} \cite{Pokhariyal-82} if 
\vspace{-0.4cm} 
\[
a_{0}=1,\quad a_{1}=-\,a_{3}=-\,\frac{1}{(n-1)},\quad
a_{2}=a_{4}=a_{5}=a_{6}=a_{7}=0, 
\]%
\vspace{-0.5cm}

\item the ${\cal W}_{9}${\em -curvature tensor} \cite{Pokhariyal-82} if 
\vspace{-0.4cm} 
\[
a_{0}=1,\quad a_{3}=-\,a_{4}=\frac{1}{(n-1)},\quad
a_{1}=a_{2}=a_{5}=a_{6}=a_{7}=0. 
\]
\end{enumerate}

Denoting 
\[
{\cal T}\left( X,Y,Z,V\right) =g({\cal T}\left( X,Y\right) Z,V), 
\]%
we write the curvature tensor ${\cal T}$ in its $\left( 0,4\right) $ form as
follows. 
\begin{eqnarray}
{\cal T}\left( X,Y,Z,V\right) &=&a_{0}\,R\left( X,Y,Z,V\right)  \nonumber \\
&&+\ a_{1}\,S\left( Y,Z\right) g\left( X,V\right) +a_{2}\,S\left( X,Z\right)
g\left( Y,V\right)  \nonumber \\
&&+\ a_{3}\,S\left( X,Y\right) g\left( Z,V\right) +a_{4}\,S\left( X,V\right)
g\left( Y,Z\right)  \nonumber \\
&&+\ a_{5}\,S\left( Y,V\right) g\left( X,Z\right) +a_{6}\,S\left( Z,V\right)
g\left( X,Y\right)  \nonumber \\
&&+\ a_{7}\,r\left( g\left( Y,Z\right) g\left( X,V\right) -g\left(
X,Z\right) g\left( Y,V\right) \right) .  \label{eq-gen-cur-1}
\end{eqnarray}
In a semi-Riemannian manifold $\left( M,g\right) $, let $\left\{
e_{i}\right\} $, $i=1,\ldots ,n$ be a local orthonormal basis, define 
\[
\left( {\rm div}\,{\cal T}\right) (X,Y,Z)=\sum_{i=1}^{n}\varepsilon
_{i}g((\nabla _{e_{i}}{\cal T)}(X,Y)Z,e_{i}), 
\]
where $\varepsilon _{i}=g(e_{i},e_{i})$. Then 
\begin{eqnarray}
\left( {\rm div}\,{\cal T}\right) (X,Y,Z) &=&(a_{0}+a_{1})(\nabla
_{X}S)(Y,Z)+(-a_{0}+a_{2})(\nabla _{Y}S)(X,Z)  \nonumber \\
&&+\,a_{3}(\nabla _{Z}S)(X,Y)+\left( \frac{a_{4}}{2}+a_{7}\right) (\nabla
_{X}r)g(Y,Z)  \nonumber \\
&&+\left( \frac{a_{5}}{2}-a_{7}\right) (\nabla _{Y}r)g(X,Z)+\frac{a_{6}}{2}%
(\nabla _{Z}r)g(X,Y),  \label{eq-divT-11}
\end{eqnarray}

\begin{eqnarray}
S_{{\cal T}}(X,Y) &=&(a_{0}+na_{1}+a_{2}+a_{3}+a_{5}+a_{6})S(X,Y)  \nonumber
\\
&&+\ (a_{4}+(n-1)a_{7})r\,g(X,Y).  \label{eq-ric-T}
\end{eqnarray}

\begin{defn-new}
An $n$-dimensional semi-Riemannian manifold is said to be ${\cal T}$%
-conservative \cite{Tripathi-Gupta} if ${\rm div}\,{\cal T}=0$.
\end{defn-new}

\begin{notation}
We will call ${\cal T}$-curvature tensor as ${\cal T}_{a}$-curvature tensor,
whenever it is necessary. If $a_{0},\ldots ,a_{7}$ are replaced by $%
b_{0},\ldots ,b_{7}$ in the definition of ${\cal T}$-curvature tensor, then
we will call ${\cal T}$-curvature tensor as ${\cal T}_{b}$-curvature tensor.
\end{notation}

\section{$(N(k),\protect\xi)$-semi-Riemannian manifolds\label{sect-NK}}

Let $(M,g)$ be an $n$-dimensional semi-Riemannian manifold \cite{ONeill-83}
equipped with a semi-Riemannian metric $g$. If ${\rm index}(g)=1$ then $g$
is a Lorentzian metric and $(M,g)$ a Lorentzian manifold \cite%
{Beem-Ehrlich-81}. If $g$ is positive definite then $g$ is an usual
Riemannian metric and $(M,g)$ a Riemannian manifold. \medskip

The $k${\em -nullity distribution} \cite{Tanno-88} of $(M,g)$ for a real
number $k$ is the distribution%
\[
N(k):p\mapsto N_{p}(k)=\left\{ Z\in
T_{p}M:R(X,Y)Z=k(g(Y,Z)X-g(X,Z)Y)\right\} . 
\]%
Let $\xi $ be a non-null unit vector field in $(M,g)$ and $\eta $ its
associated $1$-form. Thus 
\[
g(\xi ,\xi )=\varepsilon , 
\]%
where $\varepsilon =1$ or $-\,1$ according as $\xi $ is spacelike or
timelike, and 
\begin{equation}
\eta \left( X\right) =\varepsilon g\left( X,\xi \right) ,\qquad \eta \left(
\xi \right) =1.  \label{eq-cond}
\end{equation}

\begin{defn-new}
An $(N(k),\xi )${\em -semi-Riemannian manifold} consists of a
semi-Riemannian manifold $(M,g)$, a $k$-nullity distribution $N(k)$ on $%
(M,g) $ and a non-null unit vector field $\xi $ in $(M,g)$ belonging to $%
N(k) $.
\end{defn-new}

Now, we intend to give some examples of $(N(k),\xi )$-semi-Riemannian
manifolds. For this purpose we collect some definitions from the geometry of
almost contact manifolds and almost paracontact manifolds as follows:

\subsection*{Almost contact manifolds}

Let $M$ be a smooth manifold of dimension $n=2m+1$. Let $\varphi $, $\xi $
and $\eta $ be tensor fields of type $(1,1)$, $(1,0)$ and $(0,1)$,
respectively. If $\varphi $, $\xi $ and $\eta $ satisfy the conditions 
\begin{equation}
\varphi ^{2}=-I+\eta \otimes \xi ,  \label{eq-str-1}
\end{equation}%
\begin{equation}
\eta (\xi )=1,  \label{eq-str-2}
\end{equation}%
where $I$ denotes the identity transformation, then $M$ is said to have an
almost contact structure $\left( \varphi ,\xi ,\eta \right) $. A manifold $M$
alongwith an almost contact structure is called an {\em almost contact
manifold} \cite{Blair-76}. Let $g$ be a semi-Riemannian metric on $M$ such
that 
\begin{equation}
g\left( \varphi X,\varphi Y\right) =g\left( X,Y\right) -\varepsilon \eta
(X)\eta \left( Y\right) ,  \label{eq-str-3}
\end{equation}%
where $\varepsilon =\pm 1$. Then $(M,g)$ is an $\left( \varepsilon \right) $-%
{\em almost contact metric manifold} \cite{Duggal-90-IJMMS} equipped with an 
$\left( \varepsilon \right) ${\em -almost contact metric structure} $%
(\varphi ,\xi ,\eta ,g,\varepsilon )$. In particular, if the metric $g$ is
positive definite, then an $(\varepsilon )$-almost contact metric manifold
is the usual {\em almost contact metric manifold }\cite{Blair-76}. \medskip

From (\ref{eq-str-3}), it follows that 
\begin{equation}
g\left( X,\varphi Y\right) =-g\left( \varphi X,Y\right)  \label{eq-str-5}
\end{equation}%
and 
\begin{equation}
g\left( X,\xi \right) =\varepsilon \eta (X).  \label{eq-str-6}
\end{equation}%
From (\ref{eq-str-2}) and (\ref{eq-str-6}), we have 
\begin{equation}
g\left( \xi ,\xi \right) =\varepsilon .  \label{eq-str-7}
\end{equation}

In an $\left( \varepsilon \right) $-almost contact metric manifold, the
fundamental $2$-form $\Phi $ is defined by 
\begin{equation}
\Phi (X,Y)=g(X,\varphi Y).  \label{eq-str-4}
\end{equation}%
An $\left( \varepsilon \right) $-almost contact metric manifold with $\Phi
=d\eta $ is an $\left( \varepsilon \right) ${\em -contact metric manifold }%
\cite{Takahashi-69}. For $\varepsilon =1$ and $g$ Riemannian, $M$ is the
usual{\em \ contact metric manifold }\cite{Blair-76}. A contact metric
manifold with $\xi \in N(k)$, is called a {\em $N(k)$-contact metric
manifold }\cite{Blair-Kim-Tripathi-05}. \medskip

An $\left( \varepsilon \right) $-almost contact metric structure $(\varphi
,\xi ,\eta ,g,\varepsilon )$ is called an $\left( \varepsilon \right) ${\em %
-Sasakian structure} if 
\[
\left( \nabla _{X}\varphi \right) Y=g(X,Y)\xi -\varepsilon \eta \left(
Y\right) X, 
\]%
where $\nabla $ is Levi-Civita connection with respect to the metric $g$. A
manifold endowed with an $\left( \varepsilon \right) $-Sasakian structure is
called an $\left( \varepsilon \right) ${\em -Sasakian manifold }\cite%
{Takahashi-69}. For $\varepsilon =1$ and $g$ Riemannian, $M$ is the usual%
{\em \ Sasakian manifold }\cite{Sasaki-60,Blair-76}. \medskip

An almost contact metric manifold is a {\em Kenmotsu manifold} \cite%
{Kenmotsu-72} if 
\begin{equation}
\left( \nabla _{X}\varphi \right) Y=g(\varphi X,Y)\xi -\eta \left( Y\right)
\varphi X.  \label{eq-str-8}
\end{equation}%
By (\ref{eq-str-8}), we have 
\begin{equation}
\nabla _{X}\xi =X-\eta (X)\xi .  \label{eq-str-9}
\end{equation}

\subsection*{Almost paracontact manifolds}

Let $M$ be an $n$-dimensional {\em almost paracontact manifold} \cite%
{Sato-76} equipped with an {\em almost paracontact structure} $\left(
\varphi ,\xi ,\eta \right) $, where $\varphi $, $\xi $ and $\eta $ are
tensor fields of type $(1,1)$, $(1,0)$ and $(0,1)$, respectively; and
satisfy the conditions 
\begin{equation}
\varphi ^{2}=I-\eta \otimes \xi ,  \label{eq-str-11}
\end{equation}%
\begin{equation}
\eta (\xi )=1.  \label{eq-str-12}
\end{equation}%
Let $g$ be a semi-Riemannian metric on $M$ such that 
\begin{equation}
g\left( \varphi X,\varphi Y\right) =g\left( X,Y\right) -\varepsilon \eta
(X)\eta \left( Y\right) ,  \label{eq-str-13}
\end{equation}%
where $\varepsilon =\pm 1$. Then $\left( M,g\right) $ is an $\left(
\varepsilon \right) ${\em -almost paracontact metric manifold} equipped with
an $\left( \varepsilon \right) ${\em -almost paracontact metric structure} $%
(\varphi ,\xi ,\eta ,g,\varepsilon )$. In particular, if ${\rm index}(g)=1$,
then an $(\varepsilon )$-almost paracontact metric manifold is said to be a 
{\em Lorentzian almost paracontact manifold}. In particular, if the metric $%
g $ is positive definite, then an $(\varepsilon )$-almost paracontact metric
manifold is the usual {\em almost paracontact metric manifold} \cite{Sato-76}%
. \medskip

The equation (\ref{eq-str-13}) is equivalent to 
\begin{equation}
g\left( X,\varphi Y\right) =g\left( \varphi X,Y\right)  \label{eq-str-14}
\end{equation}%
along with 
\begin{equation}
g\left( X,\xi \right) =\varepsilon \eta (X).  \label{eq-str-15}
\end{equation}%
From (\ref{eq-str-12}) and (\ref{eq-str-15}), we have 
\begin{equation}
g\left( \xi ,\xi \right) =\varepsilon .  \label{eq-str-16}
\end{equation}

An $\left( \varepsilon \right) $-almost paracontact metric structure is
called an $\left( \varepsilon \right) ${\em -para-Sasakian structure} \cite%
{TKYK-09} if 
\begin{equation}
\left( \nabla _{X}\varphi \right) Y=-\,g(\varphi X,\varphi Y)\xi
-\varepsilon \eta \left( Y\right) \varphi ^{2}X,  \label{eq-str-17}
\end{equation}%
where $\nabla $ is Levi-Civita connection with respect to the metric $g$. A
manifold endowed with an $\left( \varepsilon \right) $-para-Sasakian
structure is called an $\left( \varepsilon \right) ${\em -para-Sasakian
manifold} \cite{TKYK-09}. For $\varepsilon =1$ and $g$ Riemannian, $M$ is
the usual para-Sasakian manifold \cite{Sato-76}. For $\varepsilon =-1$, $g$
Lorentzian and $\xi $ replaced by $-\xi $, $M$ becomes a Lorentzian
para-Sasakian manifold \cite{Matsumoto-89}. \bigskip

\begin{ex-new}
The following are some well known examples of $(N(k),\xi )$-semi-Riemannian
manifolds:

\begin{enumerate}
\item An $N(k)$-contact metric manifold \cite{Blair-Kim-Tripathi-05} is an $%
(N(k),\xi )$-Riemannian manifold.

\item A Sasakian manifold \cite{Sasaki-60} is an $(N(1),\xi )$-Riemannian
manifold.

\item A Kenmotsu manifold \cite{Kenmotsu-72} is an $(N(-1),\xi )$-Riemannian
manifold.

\item An $(\varepsilon )$-Sasakian manifold \cite{Takahashi-69} an $%
(N(\varepsilon ),\xi )$-semi-Riemannian manifold.

\item A para-Sasakian manifold \cite{Sato-76} is an $(N(-1),\xi )$%
-Riemannian manifold.

\item An $(\varepsilon )$-para-Sasakian manifold \cite{TKYK-09} is an $%
(N(-\varepsilon ),\xi )$-semi-Riemannian manifold.
\end{enumerate}
\end{ex-new}

In an $n$-dimensional $\left( N(k),\xi \right) $-semi-Riemannian manifold $%
(M,g)$, it is easy to verify that 
\begin{equation}
R(X,Y)\xi =\varepsilon k(\eta (Y)X-\eta (X)Y),  \label{eq-curvature}
\end{equation}%
\begin{equation}
R(\xi ,X)Y=\varepsilon k(\varepsilon g(X,Y)\xi -\eta (Y)X),
\label{eq-curvature-2}
\end{equation}%
\begin{equation}
R(\xi ,X)\xi =\varepsilon k(\eta (X)\xi -X),  \label{eq-curvature-3}
\end{equation}%
\begin{equation}
R\left( X,Y,Z,\xi \right) =\varepsilon k(\,\eta \left( X\right) g\left(
Y,Z\right) -\eta \left( Y\right) g\left( X,Z\right) ),
\label{eq-eps-PS-R(X,Y,Z,xi)}
\end{equation}%
\begin{equation}
\eta \left( R\left( X,Y\right) Z\right) =k(\eta \left( X\right) g\left(
Y,Z\right) -\eta \left( Y\right) g\left( X,Z\right) ),
\label{eq-eps-PS-eta(R(X,Y),Z)}
\end{equation}%
\begin{equation}
S(X,\xi )=\varepsilon k(n-1)\eta (X),  \label{eq-ricci}
\end{equation}%
\begin{equation}
Q\xi =k(n-1)\xi ,  \label{eq-Q}
\end{equation}%
\begin{equation}
S(\xi ,\xi )=\varepsilon k(n-1),  \label{eq-S-xi-xi}
\end{equation}%
\begin{equation}
\eta (QX)=\varepsilon g(QX,\xi )=\varepsilon S(X,\xi )=k(n-1)\eta (X).
\label{eq-eta-QX}
\end{equation}%
Moreover, define 
\begin{equation}
S^{\ell }(X,Y)=g(Q^{\ell }X,Y)=S(Q^{\ell -1}X,Y),  \label{eq-S-p}
\end{equation}%
where $\ell =0,1,2,\ldots $ and $S^{0}=g$. Using (\ref{eq-eta-QX}) in (\ref%
{eq-S-p}), we get 
\begin{equation}
S^{\ell }(X,\xi )=\varepsilon k^{\ell }(n-1)^{\ell }\eta (X).
\label{eq-Sp-QX-xi}
\end{equation}

Now, we state the following Lemma without proof.

\begin{lem}
\label{GCT} Let $M$ be an $n$-dimensional $\left( N(k),\xi \right) $%
-semi-Riemannian manifold. Then 
\begin{eqnarray}
{\cal T}_{a}(X,Y)\xi &=&(-\varepsilon ka_{0}+\varepsilon
k(n-1)a_{2}-\varepsilon a_{7}\,r)\eta (X)Y  \nonumber \\
&&+\,(\varepsilon ka_{0}+\varepsilon k(n-1)a_{1}+\varepsilon a_{7}\,r)\,\eta
(Y)X  \nonumber \\
&&+\,a_{3}\,S(X,Y)\xi +\varepsilon a_{4}\,\eta (Y)QX  \nonumber \\
&&+\,\varepsilon a_{5}\eta (X)QY+k(n-1)a_{6}g(X,Y)\xi ,  \label{eq-X-Y-xi}
\end{eqnarray}%
\begin{eqnarray}
{\cal T}_{\!a}(\xi ,X)\xi &=&(-\varepsilon ka_{0}\,+\,\varepsilon
k(n-1)a_{2}-\varepsilon a_{7}\,r)X+\varepsilon a_{5}\,QX  \nonumber \\
&&+\left\{ \varepsilon ka_{0}+\varepsilon k(n-1)a_{1}+\varepsilon
k(n-1)a_{3}\right.  \nonumber \\
&&\quad \quad +\left. \varepsilon k(n-1)a_{4}+\varepsilon
k(n-1)a_{6}\,+\varepsilon a_{7}\,r\right\} \eta (X)\xi ,  \label{eq-xi-X-xi}
\end{eqnarray}%
\begin{eqnarray}
{\cal T}_{\!a}(\xi ,Y)Z &=&(ka_{0}+k(n-1)a_{4}+a_{7}r)g\left( Y,Z\right) \xi
\,  \nonumber \\
&&+\,a_{1}\,S\left( Y,Z\right) \xi +\varepsilon k(n-1)a_{3}\eta (Y)Z 
\nonumber \\
&&+\,\varepsilon a_{5}\,\eta (Z)QY+\varepsilon a_{6}\,\eta (Y)QZ  \nonumber
\\
&&+\,(-\varepsilon ka_{0}+\varepsilon k(n-1)a_{2}\,-\varepsilon
a_{7}\,r)\eta (Z)Y,  \label{eq-xi-Y-Z}
\end{eqnarray}%
\begin{eqnarray}
\eta ({\cal T}_{\!a}\left( X,Y)\xi \right) &=&\varepsilon
k(n-1)(a_{1}+a_{2}+a_{4}+a_{5})\eta (X)\eta (Y)  \nonumber \\
&&+\,a_{3}\,S(X,Y)+\,k(n-1)a_{6}g(X,Y),  \label{eq-eta-xi-X-Y}
\end{eqnarray}%
\begin{eqnarray}
{\cal T}_{\!a}\left( X,Y,\xi ,V\right) &=&(-\varepsilon ka_{0}\,+\varepsilon
k(n-1)a_{2}\,-\varepsilon a_{7}\,r)\eta (X)g(Y,V)  \nonumber \\
&&+\ (\varepsilon ka_{0}+\varepsilon k(n-1)a_{1}+\varepsilon a_{7}\,r)\,\eta
(Y)g(X,V)  \nonumber \\
&&+\,\varepsilon a_{3}\,S(X,Y)\eta (V)+\varepsilon a_{4}\,\eta (Y)S(X,V) 
\nonumber \\
&&+\,\varepsilon a_{5}\,\eta (X)S(Y,V)+\varepsilon k(n-1)a_{6}\,g(X,Y)\eta
(V),  \label{eq-X-Y-xi-V}
\end{eqnarray}%
\begin{eqnarray}
{\cal T}_{\!a}(X,\xi )\xi &=&\left\{ -\varepsilon ka_{0}\,+\varepsilon
k(n-1)a_{2}\,+\varepsilon k(n-1)a_{3}\right.  \nonumber \\
&&\quad +\left. \varepsilon k(n-1)a_{5}+\varepsilon
k(n-1)a_{6}\,\,-\varepsilon a_{7}\,r\right\} \eta (X)\xi  \nonumber \\
&&+\ (\varepsilon ka_{0}+\varepsilon k(n-1)a_{1}+\varepsilon
a_{7}\,r)\,X+\varepsilon a_{4}\,QX,  \label{eq-X-xi-xi}
\end{eqnarray}%
\begin{eqnarray}
S_{{\cal T}_{\!a}}(X,\xi ) &=&\left\{ \varepsilon
k(n-1)(a_{0}+na_{1}+a_{2}+a_{3}+a_{5}+a_{6})\right.  \nonumber \\
&&\left. +\,\varepsilon r(a_{4}+(n-1)a_{7})\right\} \eta (X),
\label{eq-ric-T1}
\end{eqnarray}%
\begin{eqnarray}
S_{{\cal T}_{\!a}}(\xi ,\xi ) &=&\varepsilon
k(n-1)(a_{0}+na_{1}+a_{2}+a_{3}+a_{5}+a_{6})  \nonumber \\
&&+\,\varepsilon r(a_{4}+(n-1)a_{7}).  \label{eq-ric-T2}
\end{eqnarray}
\end{lem}

\begin{rem-new}
The relations (\ref{eq-curvature}) -- (\ref{eq-ric-T2}) are true for

\begin{enumerate}
\item a $N(k)$-contact metric manifold \cite{Blair-Kim-Tripathi-05}\ ($%
\varepsilon =1$),

\item a Sasakian manifold \cite{Sasaki-60} ($k=1$, $\varepsilon =1$),

\item a Kenmotsu manifold \cite{Kenmotsu-72} ($k=-1$, $\varepsilon =1$),

\item an $(\varepsilon )$-Sasakian manifold \cite{Takahashi-69} ($%
k=\varepsilon $, $\varepsilon k=1$),

\item a para-Sasakian manifold \cite{Sato-76} ($k=-1$, $\varepsilon =1$), and

\item an $(\varepsilon )$-para-Sasakian manifold \cite{TKYK-09} ($%
k=-\,\varepsilon $, $\varepsilon k=-\,1$).
\end{enumerate}

\noindent Even, all the relations and results of this paper will be true for
the above six cases.
\end{rem-new}

\section{{$\protect\xi $-${\cal T}_{\!a}$-flat $(N(k),\protect\xi)$%
-semi-Riemannian manifolds \label{sect-T0}}}

\begin{defn-new}
An $n$-dimensional $(N(k),\xi )$-semi-Riemannian manifold $(M,g)$ is said to
be $\xi $-${\cal T}_{\!a}$-flat if it satisfies 
\[
{\cal T}_{\!a}(X,Y)\xi =0.
\]%
In particular, if {\em ${\cal T}_{\!a}$} is equal to $R$, ${\cal C}_{\ast }$%
, ${\cal C}$, ${\cal L}$, ${\cal V}$, ${\cal P}_{\ast }$, ${\cal P}$, ${\cal %
M}$, ${\cal W}_{0}$, ${\cal W}_{0}^{\ast }$, ${\cal W}_{1}$, ${\cal W}%
_{1}^{\ast }$, ${\cal W}_{2}$, ${\cal W}_{3}$, ${\cal W}_{4}$, ${\cal W}_{5}$%
, ${\cal W}_{6}$, ${\cal W}_{7}$, ${\cal W}_{8}$, ${\cal W}_{9}$, then it
becomes $\xi $-flat, $\xi $-quasi-conformally flat, $\xi $-conformally{\bf \ 
}flat, $\xi $-conharmonically flat, $\xi $-concircularly flat, $\xi $%
-pseudo-projectively flat, $\xi $-projectively flat, $\xi $-${\cal M}$-flat, 
$\xi $-${\cal W}_{0}$-flat, $\xi $-${\cal W}_{0}^{\ast }$-flat, $\xi $-$%
{\cal W}_{1}$-flat, $\xi $-${\cal W}_{1}^{\ast }$-flat, $\xi $-${\cal W}_{2}$%
-flat, $\xi $-${\cal W}_{3}$-flat, $\xi $-${\cal W}_{4}$-flat, $\xi $-${\cal %
W}_{5}$-flat, $\xi $-${\cal W}_{6}$-flat, $\xi $-${\cal W}_{7}$-flat, $\xi $-%
${\cal W}_{8}$-flat, $\xi $-${\cal W}_{9}$-flat, respectively.
\end{defn-new}

\begin{th}
\label{th-11} Let $M$ be an $n$-dimensional $\xi $-${\cal T}_{\!a}$-flat $%
(N(k),\xi )$-semi-Riemannian manifold.

\begin{enumerate}
\item[{\rm 1.}] If $a_{4}\not=0$ and $a_{4}+(n-1)a_{7}\neq 0$, then 
\begin{equation}
S=D_{2}\,g+D_{3}\eta \otimes \eta ,  \label{eq-xi-Tflat-6-11}
\end{equation}%
where 
\begin{equation}
D_{2}=-\,\frac{ka_{0}+k(n-1)a_{1}+a_{7}r}{a_{4}}  \label{eq-xi-Tflat-6i-11}
\end{equation}%
and 
\begin{equation}
D_{3}=\frac{ka_{0}-(n-1)k(a_{2}+a_{3}+a_{5}+a_{6})+a_{7}r}{\varepsilon a_{4}}%
.  \label{eq-xi-Tflat-6ii-11}
\end{equation}

Therefore it is an $\eta $-Einstein manifold and 
\begin{equation}
r=-\,\frac{k(n-1)\left( a_{0}+na_{1}+a_{2}+a_{3}+a_{5}+a_{6}\right) }{%
a_{4}+(n-1)a_{7}}=D_{1} \,{\rm (say)}.  \label{eq-r-11}
\end{equation}
In particular, $M$ becomes an Einstein manifold provided 
\[
ka_{0}-(n-1)k(a_{2}+a_{3}+a_{5}+a_{6})+a_{7}r=0. 
\]

\item[{\rm 2.}] If $a_{4}=0$ and $a_{7}\not=0$, then 
\begin{equation}
r=-\,\frac{k\left( a_{0}+na_{1}+a_{2}+a_{3}+a_{5}+a_{6}\right) }{a_{7}}.
\label{eq-T-r11-11}
\end{equation}

\item[{\rm 3.}] If $a_{4}=0$ and $a_{7}=0$, then either $k=0$ or 
\begin{equation}
a_{0}+na_{1}+a_{2}+a_{3}+a_{5}+a_{6}=0.  \label{eq-T-r12-11}
\end{equation}
\end{enumerate}
\end{th}

\noindent {\bf Proof.} By (\ref{eq-gen-cur-1}) and (\ref{eq-cond}), we get 
\begin{eqnarray}
{\cal T}_{\!a}(X,Y,\xi ,W) &=&a_{0}\,R(X,Y,\xi ,W)  \nonumber \\
&&+\ a_{1}\,S\left( Y,\xi \right) g(X,W)+a_{2}S\left( X,\xi \right) g(Y,W) 
\nonumber \\
&&+\ a_{3}\varepsilon S(X,Y)\eta (W)+a_{4}\varepsilon \eta (Y)S(X,W) 
\nonumber \\
&&+\ a_{5}\varepsilon \eta (X)S(Y,W)+a_{6}\,g(X,Y)S(\xi ,W)  \nonumber \\
&&+\ a_{7}\varepsilon \,r\left( \eta (Y)g(X,W)-\eta (X)g(Y,W)\right) .
\label{eq-xi-Tflat-1-11}
\end{eqnarray}%
Using $Y=\xi $ in (\ref{eq-xi-Tflat-1-11}), we get 
\begin{eqnarray}
{\cal T}_{\!a}(X,\xi ,\xi ,W) &=&a_{0}\,R(X,\xi ,\xi ,W)  \nonumber \\
&&+\,a_{1}\,S\left( \xi ,\xi \right) g(X,W)+a_{2}\,S\left( X,\xi \right)
g(\xi ,W)  \nonumber \\
&&+\,a_{3}\varepsilon S(X,\xi )\eta (W)+\ a_{4}\varepsilon \eta (\xi )S(X,W)
\nonumber \\
&&+\,a_{5}\varepsilon \eta (X)S(\xi ,W)+a_{6}\,g(X,\xi )S(\xi ,W)  \nonumber
\\
&&+\,a_{7}\varepsilon r\left( \eta (\xi )g(X,W)-\eta (X)g(\xi ,W)\right) .
\label{eq-xi-Tflat-5-11}
\end{eqnarray}%
{\bf Case 1.} If $a_{4}\not=0$ and $a_{4}+(n-1)a_{7}\neq 0$, then using (\ref%
{eq-cond}), (\ref{eq-curvature-3}) (\ref{eq-ricci}) and (\ref{eq-S-xi-xi})
and the fact that $M$ is $\xi $-${\cal T}_{\!a}$-flat in (\ref%
{eq-xi-Tflat-5-11}), we get (\ref{eq-xi-Tflat-6-11}) and (\ref{eq-r-11}).

\noindent {\bf Case 2.} If $a_{4}=0$ and $a_{7}\not=0$, then by using (\ref%
{eq-cond}), (\ref{eq-curvature-3}) (\ref{eq-ricci}) and (\ref{eq-S-xi-xi})
and the fact that $M$ is $\xi $-${\cal T}_{\!a}$-flat in (\ref%
{eq-xi-Tflat-5-11}), we get 
\begin{eqnarray*}
\varepsilon \left( a_{0}k+(n-1)ka_{1}+a_{7}r\right) g(Y,Z) &=&\left(
-(n-1)k(a_{2}+a_{3}+a_{5}+a_{6})\right. \\
&&+\,\left. (a_{0}k+a_{7}\,r)\right) \eta (Y)\eta (Z).
\end{eqnarray*}%
Contracting the above equation, we get 
\begin{equation}
a_{7}r=-k\left( a_{0}+na_{1}+a_{2}+a_{3}+a_{5}+a_{6}\right) .
\label{eq-T-r13-11}
\end{equation}%
Since $a_{7}\not=0$, we get (\ref{eq-T-r11-11}).

\noindent {\bf Case 3.} If $a_{4}=0$ and $a_{7}=0$, then from (\ref%
{eq-T-r13-11}) either $k=0$ or (\ref{eq-T-r12-11}) is satisfied. This proves
the result. $\blacksquare $

\begin{th}
\label{th-xi} Let $M$ be an $n$-dimensional $(N(k),\xi )$-semi-Riemannian
manifold such that $a_{4}\not=0$ and $a_{4}+(n-1)a_{7}\not=0$. If $M$
satisfies $(\ref{eq-xi-Tflat-6-11})$, then 
\begin{eqnarray}
{\cal T}_{\!a}(X,Y)\xi &=&(\varepsilon ka_{0}+\varepsilon
k(n-1)a_{1}+\varepsilon a_{4}D_{2}+\varepsilon a_{7}D_{1})\eta (Y)X 
\nonumber \\
&&+\,(-\varepsilon ka_{0}+\varepsilon k(n-1)a_{2}+\varepsilon
a_{5}D_{2}-\varepsilon a_{7}D_{1})\eta (X)Y  \nonumber \\
&&+\,(k(n-1)a_{6}+a_{3}D_{2})g(X,Y)\xi  \nonumber \\
&&+\,(a_{3}+a_{4}+a_{5})D_{3}\eta (X)\eta (Y)\xi .  \label{eq-RRR}
\end{eqnarray}
\end{th}

\begin{rem-new}
If $M$ is $\xi $-conformally flat $(N(k),\xi )$-semi-Riemannian manifold,
then from {\rm (\ref{eq-r-11})}, the scalar curvature $r$ is in
indeterminate form.
\end{rem-new}

Suppose that a $\left( N(k),\xi \right) $-semi-Riemannian manifold is $\eta $%
-Einstein. Then there are functions $\alpha $ and $\beta $\ such that 
\begin{equation}
S(X,Y)=\alpha g(X,Y)+\beta \eta (X)\eta (Y).  \label{eq-T-flat-S1}
\end{equation}%
On contracting (\ref{eq-T-flat-S1}), we get 
\begin{equation}
r=\alpha n+\beta \varepsilon .  \label{eq-1}
\end{equation}%
Taking $X=\xi =Y$ in (\ref{eq-T-flat-S1}), we get 
\begin{equation}
k(n-1)=\alpha +\beta \varepsilon .  \label{eq-2}
\end{equation}%
Using (\ref{eq-1}) in (\ref{eq-2}) yields 
\begin{equation}
r=(k+\alpha )(n-1).  \label{eq-3}
\end{equation}

\begin{th}
Let $M$ be an $\eta $-Einstein $\left( N(k),\xi \right) $-semi-Riemannian
manifold. Then 
\begin{eqnarray}
{\cal T}_{\!a}(X,Y,\xi ,V) &=&\varepsilon \left( ka_{0}+(\alpha +\beta
)a_{1}+\alpha a_{4}+(k+\alpha )(n-1)a_{7}\right) \eta (Y)g(X,V)  \nonumber \\
&&+\,\varepsilon \left( -ka_{0}+(\alpha +\beta )a_{2}+\alpha a_{5}-(k+\alpha
)(n-1)a_{7}\right) \eta (X)g(Y,V)  \nonumber \\
&&+\,\varepsilon \left( \alpha a_{3}+(\alpha +\beta )a_{6}\right) \eta
(V)g(X,Y)  \nonumber \\
&&+\,\varepsilon \beta (a_{3}+a_{4}+a_{5})\eta (X)\eta (Y)\eta (V).
\label{eq-T-XYxi}
\end{eqnarray}
\end{th}

\noindent {\bf Proof.} Let $M$ be an $\eta $-Einstein $\left( N(k),\xi
\right) $-semi-Riemannian manifold. Taking $Z=\xi $ in (\ref{eq-gen-cur-1})
and the using (\ref{eq-curvature}), (\ref{eq-ricci}), (\ref{eq-Q}) and (\ref%
{eq-3}), we get (\ref{eq-T-XYxi}). $\blacksquare $

In view of Theorem~\ref{th-11}, we have the following Corollaries:

\begin{cor}
Let $M$ be an $n$-dimensional $\xi $-quasi-conformally flat $\left( N(k),\xi
\right) $-semi-Riemannian manifold such that $a_{1}\not=0$ and $%
a_{0}+(n-2)a_{1}\not=0$. Then we have the following table\/{\rm :}~%
\[
\begin{tabular}{|l|l|}
\hline
${\boldmath M}$ & ${\boldmath S=}$ \\ \hline
$N(k)$-contact metric & $k(n-1)g$ \\ \hline
Sasakian & $(n-1)g$ \\ \hline
Kenmotsu & $-\,(n-1)g$ \\ \hline
$(\varepsilon )$-Sasakian & $\varepsilon (n-1)g$ \\ \hline
para-Sasakian & $-\,(n-1)g$ \\ \hline
$(\varepsilon )$-para-Sasakian & $-\,\varepsilon (n-1)g$ \\ \hline
\end{tabular}%
\]
\end{cor}

\begin{cor}
\label{cor-c} Let $M$ be an $n$-dimensional $\xi $-conformally flat $%
(N(k),\xi )$-semi-Riemannian manifold. Then we have the following table\/%
{\rm :}~ 
\[
\begin{tabular}{|l|l|}
\hline
${\boldmath M}$ & ${\boldmath S=}$ \\ \hline
$N(k)$-contact metric {\rm \cite{De-Biswas-06}} & $\left( \dfrac{r}{n-1}%
-k\right) g+\left( nk-\dfrac{r}{n-1}\right) \eta \otimes \eta $ \\ \hline
Sasakian & $\left( \dfrac{r}{n-1}-1\right) g+\left( n-\dfrac{r}{n-1}\right)
\eta \otimes \eta $ \\ \hline
Kenmotsu & $\left( \dfrac{r}{n-1}+1\right) g-\left( n+\dfrac{r}{n-1}\right)
\eta \otimes \eta $ \\ \hline
$(\varepsilon )$-Sasakian & $\left( \dfrac{r}{n-1}-\varepsilon \right)
g+\varepsilon \left( \varepsilon n-\dfrac{r}{n-1}\right) \eta \otimes \eta $
\\ \hline
para-Sasakian & $\left( \dfrac{r}{n-1}+1\right) g-\left( n+\dfrac{r}{n-1}%
\right) \eta \otimes \eta $ \\ \hline
$(\varepsilon )$-para-Sasakian & $\left( \dfrac{r}{n-1}+\varepsilon \right)
g-\varepsilon \left( \varepsilon n+\dfrac{r}{n-1}\right) \eta \otimes \eta $
\\ \hline
\end{tabular}%
\ 
\]
\end{cor}

\begin{cor}
Let $M$ be an $n$-dimensional $\xi $-conharmonically flat $(N(k),\xi )$%
-semi-Riemannian manifold. Then we have the following table\/{\rm :}~%
\[
\begin{tabular}{|l|l|}
\hline
${\boldmath M}$ & ${\boldmath S=}$ \\ \hline
$N(k)$-contact metric & $-\,kg+kn\eta \otimes \eta $ \\ \hline
Sasakian {\rm \cite{Dwivedi-Kim-10}} & $-\,g+n\eta \otimes \eta $ \\ \hline
Kenmotsu & $g-n\eta \otimes \eta $ \\ \hline
$(\varepsilon )$-Sasakian & $-\,\varepsilon g+n\eta \otimes \eta $ \\ \hline
para-Sasakian & $g-n\eta \otimes \eta $ \\ \hline
$(\varepsilon )$-para-Sasakian & $\varepsilon g-n\eta \otimes \eta $ \\ 
\hline
\end{tabular}%
\ 
\]
\end{cor}

\begin{cor}
Let $M$ be an $n$-dimensional $\xi $-concircularly flat $\left( N(k),\xi
\right) $-semi-Riemannian manifold. Then we have the following table\/{\rm :}%
~%
\[
\begin{tabular}{|l|l|}
\hline
${\boldmath M}$ & ${\boldmath r=}$ \\ \hline
$N(k)$-contact metric & $kn(n-1)$ \\ \hline
Sasakian & $n(n-1)$ \\ \hline
Kenmotsu & $-\,n(n-1)$ \\ \hline
$(\varepsilon )$-Sasakian & $\varepsilon n(n-1)$ \\ \hline
para-Sasakian & $-\,n(n-1)$ \\ \hline
$(\varepsilon )$-para-Sasakian & $-\,\varepsilon n(n-1)$ \\ \hline
\end{tabular}%
\]
\end{cor}

\begin{cor}
Let $M$ be an $n$-dimensional $\xi $-pseudo-projectively flat $\left(
N(k),\xi \right) $-semi-Riemannian manifold such that $a_{0}+(n-1)a_{1}%
\not=0 $. Then we have the following table\/{\rm :}~%
\[
\begin{tabular}{|l|l|}
\hline
${\boldmath M}$ & ${\boldmath r=}$ \\ \hline
$N(k)$-contact metric & $kn(n-1)$ \\ \hline
Sasakian & $n(n-1)$ \\ \hline
Kenmotsu & $-\,n(n-1)$ \\ \hline
$(\varepsilon )$-Sasakian & $\varepsilon n(n-1)$ \\ \hline
para-Sasakian & $-\,n(n-1)$ \\ \hline
$(\varepsilon )$-para-Sasakian & $-\,\varepsilon n(n-1)$ \\ \hline
\end{tabular}%
\]
\end{cor}

\begin{cor}
Let $M$ be an $n$-dimensional $\xi $-${\cal M}$-flat $\left( N(k),\xi
\right) $-semi-Riemannian manifold. Then we have the following table\/{\rm :}%
~%
\[
\begin{tabular}{|l|l|}
\hline
${\boldmath M}$ & ${\boldmath S=}$ \\ \hline
$N(k)$-contact metric & $k(n-1)g$ \\ \hline
Sasakian & $(n-1)g$ \\ \hline
Kenmotsu & $-\,(n-1)g$ \\ \hline
$(\varepsilon )$-Sasakian & $\varepsilon (n-1)g$ \\ \hline
para-Sasakian & $-\,(n-1)g$ \\ \hline
$(\varepsilon )$-para-Sasakian & $-\,\varepsilon (n-1)g$ \\ \hline
\end{tabular}%
\]
\end{cor}

\begin{cor}
Let $M$ be an $n$-dimensional $\xi $-${\cal W}_{2}$-flat $\left( N(k),\xi
\right) $-semi-Riemannian manifold. Then we have the following table\/{\rm :}%
~%
\[
\begin{tabular}{|l|l|}
\hline
${\boldmath M}$ & ${\boldmath S=}$ \\ \hline
$N(k)$-contact metric & $k(n-1)g$ \\ \hline
Sasakian & $(n-1)g$ \\ \hline
Kenmotsu & $-\,(n-1)g$ \\ \hline
$(\varepsilon )$-Sasakian & $\varepsilon (n-1)g$ \\ \hline
para-Sasakian & $-\,(n-1)g$ \\ \hline
$(\varepsilon )$-para-Sasakian & $-\,\varepsilon (n-1)g$ \\ \hline
\end{tabular}%
\]
\end{cor}

\begin{cor}
Let $M$ be an $n$-dimensional $\xi $-${\cal W}_{3}$-flat $\left( N(k),\xi
\right) $-semi-Riemannian manifold. Then we have the following table\/{\rm :}%
~%
\[
\begin{tabular}{|l|l|}
\hline
${\boldmath M}$ & ${\boldmath S=}$ \\ \hline
$N(k)$-contact metric & $-\,k(n-1)g+2k(n-1)\eta \otimes \eta $ \\ \hline
Sasakian & $-\,(n-1)g+2(n-1)\eta \otimes \eta $ \\ \hline
Kenmotsu & $(n-1)g-2(n-1)\eta \otimes \eta $ \\ \hline
$(\varepsilon )$-Sasakian & $-\,\varepsilon (n-1)g+2(n-1)\eta \otimes \eta $
\\ \hline
para-Sasakian & $(n-1)g-2(n-1)\eta \otimes \eta $ \\ \hline
$(\varepsilon )$-para-Sasakian & $\varepsilon (n-1)g-2(n-1)\eta \otimes \eta 
$ \\ \hline
\end{tabular}%
\]
\end{cor}

\begin{cor}
Let $M$ be an $n$-dimensional $\xi $-${\cal W}_{7}$-flat $\left( N(k),\xi
\right) $-semi-Riemannian manifold. Then we have the following table\/{\rm :}%
~\ 
\[
\begin{tabular}{|l|l|}
\hline
${\boldmath M}$ & ${\boldmath S=}$ \\ \hline
$N(k)$-contact metric & $k(n-1)\eta \otimes \eta $ \\ \hline
Sasakian & $(n-1)\eta \otimes \eta $ \\ \hline
Kenmotsu & $-\,(n-1)\eta \otimes \eta $ \\ \hline
$(\varepsilon )$-Sasakian & $(n-1)\eta \otimes \eta $ \\ \hline
para-Sasakian & $-\,(n-1)\eta \otimes \eta $ \\ \hline
$(\varepsilon )$-para-Sasakian & $-\,(n-1)\eta \otimes \eta $ \\ \hline
\end{tabular}%
\ 
\]
\end{cor}

\begin{cor}
Let $M$ be an $n$-dimensional $\xi $-${\cal W}_{9}$-flat $\left( N(k),\xi
\right) $-semi-Riemannian manifold. Then we have the following table\/{\rm :}%
~ 
\[
\begin{tabular}{|l|l|}
\hline
${\boldmath M}$ & ${\boldmath S=}$ \\ \hline
$N(k)$-contact metric & $k(n-1)g$ \\ \hline
Sasakian & $(n-1)g$ \\ \hline
Kenmotsu & $-\,(n-1)g$ \\ \hline
$(\varepsilon )$-Sasakian & $\varepsilon (n-1)g$ \\ \hline
para-Sasakian & $-\,(n-1)g$ \\ \hline
$(\varepsilon )$-para-Sasakian & $-\,\varepsilon (n-1)g$ \\ \hline
\end{tabular}%
\]
\end{cor}

\begin{rem-new}
For projective curvature tensor, ${\cal W}_{0}$-curvature tensor, ${\cal W}%
_{1}^{\ast }$-curvature tensor, ${\cal W}_{6}$-curvature tensor and ${\cal W}%
_{8}$-curvature tensor, the equation $(\ref{eq-T-r12-11})$ is true. For $%
{\cal W}_{0}^{\ast }$-curvature tensor, ${\cal W}_{1}$-curvature tensor, $%
{\cal W}_{4}$-curvature tensor and ${\cal W}_{5}$-curvature tensor $(\ref%
{eq-T-r12-11})$ do not hold.
\end{rem-new}

In view of Theorem~\ref{th-11} and Theorem~\ref{th-xi} , we have the
following

\begin{cor}
Let $M$ be an $n$-dimensional $(N(K),\xi )$-semi-Riemannian manifold. Then
the following statements are true:

\begin{enumerate}
\item[{\rm (a)}] For ${\cal T}_{\!a}\in \left\{ {\cal C},{\cal L}\right\} $, 
$M$ is $\xi $-${\cal T}_{\!a}$-flat if and only if it is $\eta $-Einstein.

\item[{\rm (b)}] For ${\cal T}_{\!a}{\in }\left\{ {\cal C}_{\ast },{\cal M},%
{\cal W}_{2}\right\} $, $M$ is $\xi $-${\cal T}_{\!a}$-flat if and only if
it is Einstein.\medskip
\end{enumerate}
\end{cor}

\section{$T$-recurrent manifolds\label{sect-TR}}

\begin{defn-new}
Let $T$ be a $(1,3)$-type tensor. A semi-Riemannian manifold $(M,g)$ is said
to be $T${\em -recurrent} if it satisfies 
\begin{equation}
(\nabla _{U}T)(X,Y)Z=\alpha (U)T(X,Y)Z,  \label{eq-rec-T}
\end{equation}%
for some nonzero $1$-form $\alpha $. In particular, if $T$ is equal to $%
{\cal T}_{\!a}$, $R$, ${\cal C}_{\ast }$, ${\cal C}$, ${\cal L}$, ${\cal V}$%
, ${\cal P}_{\ast }$, ${\cal P}$, ${\cal M} $, ${\cal W}_{0}$, ${\cal W}%
_{0}^{\ast }$, ${\cal W}_{1}$, ${\cal W}_{1}^{\ast }$, ${\cal W}_{2}$, $%
{\cal W}_{3}$, ${\cal W}_{4}$, ${\cal W}_{5}$, ${\cal W}_{6}$, ${\cal W}_{7}$%
, ${\cal W}_{8}$, ${\cal W}_{9}$, then it becomes ${\cal T}_{\!a}$%
-recurrent, recurrent, quasi-conformal recurrent, Weyl recurrent,
conharmonic recurrent, concircular recurrent, pseudo-projective recurrent,
projective recurrent, ${\cal M}$-recurrent, ${\cal W}_{0}$-recurrent, ${\cal %
W}_{0}^{\ast }$-recurrent, ${\cal W}_{1}$-recurrent, ${\cal W}_{1}^{\ast }$%
-recurrent, ${\cal W}_{2}$-recurrent, ${\cal W}_{3}$-recurrent, ${\cal W}_{4}
$-recurrent, ${\cal W}_{5}$-recurrent, ${\cal W}_{6}$-recurrent, ${\cal W}%
_{7}$-recurrent, ${\cal W}_{8}$-recurrent, ${\cal W}_{9}$-recurrent,
respectively.
\end{defn-new}

\begin{defn-new}
Let $T$ be a $(1,3)$-type tensor. A semi-Riemannian manifold $(M,g)$ is said
to be $T${\em -symmetric} if it satisfies if 
\[
\nabla T=0. 
\]
In particular, if $T$ is equal to ${\cal T}_{\!a}$, $R$, ${\cal C}_{\ast }$, 
${\cal C}$, ${\cal L}$, ${\cal V}$, ${\cal P}_{\ast }$, ${\cal P}$, ${\cal M}
$, ${\cal W}_{0}$, ${\cal W}_{0}^{\ast }$, ${\cal W}_{1}$, ${\cal W}%
_{1}^{\ast }$, ${\cal W}_{2}$, ${\cal W}_{3}$, ${\cal W}_{4}$, ${\cal W}_{5}$%
, ${\cal W}_{6}$, ${\cal W}_{7}$, ${\cal W}_{8}$, ${\cal W}_{9}$, then it
becomes ${\cal T}_{\!a}$-symmetric, symmetric, quasi-conformal symmetric,
Weyl symmetric, conharmonic symmetric, concircular symmetric,
pseudo-projective symmetric, projective symmetric, ${\cal M}$-symmetric, $%
{\cal W}_{0}$-symmetric, ${\cal W}_{0}^{\ast }$-symmetric, ${\cal W}_{1}$%
-symmetric, ${\cal W}_{1}^{\ast }$-symmetric, ${\cal W}_{2}$-symmetric, $%
{\cal W}_{3}$-symmetric, ${\cal W}_{4}$-symmetric, ${\cal W}_{5}$-symmetric, 
${\cal W}_{6}$-symmetric, ${\cal W}_{7}$-symmetric, ${\cal W}_{8} $%
-symmetric, ${\cal W}_{9}$-symmetric, respectively.
\end{defn-new}

\begin{defn-new}
Let $T$ be a $(1,3)$-type tensor. A semi-Riemannian manifold $(M,g)$ is said
to be $T${\em -semisymmetric} if it satisfies if 
\[
R(V,U)\cdot T=0, 
\]%
where $R\left( V,U\right) $ acts as a derivation on $T$. In particular, if $T
$ is equal to ${\cal T}_{\!a}$, $R$, ${\cal C}_{\ast }$, ${\cal C}$, ${\cal L%
}$, ${\cal V}$, ${\cal P}_{\ast }$, ${\cal P}$, ${\cal M} $, ${\cal W}_{0}$, 
${\cal W}_{0}^{\ast }$, ${\cal W}_{1}$, ${\cal W}_{1}^{\ast }$, ${\cal W}_{2}
$, ${\cal W}_{3}$, ${\cal W}_{4}$, ${\cal W}_{5}$, ${\cal W}_{6}$, ${\cal W}%
_{7}$, ${\cal W}_{8}$, ${\cal W}_{9}$, then it becomes ${\cal T}_{\!a}$%
-semisymmetric, semisymmetric, quasi-conformal semisymmetric, Weyl
semisymmetric, conharmonic semisymmetric, concircular semisymmetric,
pseudo-projective semisymmetric, projective semisymmetric, ${\cal M}$%
-semisymmetric, ${\cal W}_{0}$-semisymmetric, ${\cal W}_{0}^{\ast }$%
-semisymmetric, ${\cal W}_{1}$-semisymmetric, ${\cal W}_{1}^{\ast }$%
-semisymmetric, ${\cal W}_{2}$-semisymmetric, ${\cal W}_{3}$-semisymmetric, $%
{\cal W}_{4}$-semisymmetric, ${\cal W}_{5}$ -semisymmetric, ${\cal W}_{6}$%
-semisymmetric, ${\cal W}_{7}$-semisymmetric, ${\cal W}_{8}$-semisymmetric, $%
{\cal W}_{9}$-semi\allowbreak symmetric, respectively.
\end{defn-new}

\begin{th}
\label{GCT-re} Let $M$ be a semi-Riemannian manifold. If $M$ is $T$%
-recurrent or $T$-symmetric then it is $T$-semisymmetric.
\end{th}

\noindent {\bf Proof.} Let us suppose that $T\neq 0$ and $M$ be a $T$%
-recurrent semi-Riemannian manifold. Then using (\ref{eq-rec-T}), we get 
\[
\nabla _{Y}(g(T,T))=2\alpha (Y)g(T,T) 
\]%
and 
\[
\nabla _{X}\nabla _{Y}(g(T,T))=2(X\alpha (Y))g(T,T)+4\alpha (X)\alpha
(Y)g(T,T), 
\]%
where the metric $g$ is extended to the inner product between the tensor
fields in the standard fashion \cite{Nomizu-Ozeki-62}. Therefore 
\begin{eqnarray*}
0 &=&(\nabla _{X}\nabla _{Y}-\nabla _{Y}\nabla _{X}-\nabla _{\lbrack
X,Y]})g(T,T) \\
&=&4d\alpha (X,Y)g(T,T).
\end{eqnarray*}%
Since $g(T,T)\neq 0$, therefore $d\alpha (X,Y)=0$. Therefore the $1$-form $%
\alpha $ is closed.

Now, from (\ref{eq-rec-T}) we have 
\[
(\nabla _{V}\nabla _{U}T)(X,Y)Z=(V\alpha (U)+\alpha (V)\alpha (U))T(X,Y)Z. 
\]%
Hence%
\[
(\nabla _{V}\nabla _{U}T-\nabla _{U}\nabla _{V}T-\nabla _{\left[ V,U\right]
}T)(X,Y)Z=\ 0. 
\]%
Therefore, we have $R(V,U)\cdot T=0$. Similarly, we can prove that if $M$ is 
$T$-symmetric semi-Riemannian manifold, then it is $T$-semisymmetric. This
proves the result. $\blacksquare $

\section{$({\cal T}_{\!a},{\cal T}_{\!b})$-semisymmetry\label{sect-TT}}

It is well known that every $(1,1)$ tensor field ${\cal A}$ on a
differentiable manifold determines a derivation ${\cal A}\cdot $ of the
tensor algebra on the manifold, commuting with contractions. For example,
the $(1,1)$ tensor fields ${\cal B}(V,U)$ induces the derivation ${\cal B}%
(V,U)\cdot $, thus associating with a $(0,s)$ tensor field ${\cal K}$, give
the $(0,s+2)$ tensor ${\cal B}\cdot {\cal K}$. The condition ${\cal B}\cdot 
{\cal K}$ is defined by%
\begin{eqnarray}
({\cal B}\cdot {\cal K})(X_{1},\ldots ,X_{s},X,Y) &=&({\cal B}(X,Y)\cdot 
{\cal K})(X_{1},\ldots ,X_{s})  \nonumber \\
&=&-\,{\cal K}({\cal B}(X,Y)X_{1},\ldots ,X_{s})-\cdots  \nonumber \\
&&-\,{\cal K}(X_{1},\ldots ,{\cal B}(X,Y)X_{s}).  \label{eq-RR}
\end{eqnarray}

\begin{defn-new}
A semi-Riemannian manifold $\left( M,g\right) $ is said to be $({\cal T}%
_{\!a},{\cal T}_{\!b})$-semisymmetric if 
\begin{equation}
{\cal T}_{\!a}(X,Y)\cdot {\cal T}_{\!b}=0,  \label{eq-T.T}
\end{equation}%
where ${\cal T}_{\!a}(X,Y)$ acts as a derivation on ${\cal T}_{\!b}$. In
particular, it is said to be $(R,{\cal T}_{\!a})$-semisymmetric if 
\begin{equation}
R(X,Y)\cdot {\cal T}_{\!a}=0,  \label{eq-R.T}
\end{equation}%
which, in brief, is said to be ${\cal T}_{\!a}$-semisymmetric.
\end{defn-new}

\begin{th}
\label{th-T-T} Let $M$ be an $n$-dimensional $({\cal T}_{\!a},{\cal T}%
_{\!b}) $-semisymmetric $(N(k),\xi )$-semi-Riemannian manifold. Then

\begin{eqnarray}
&&-\,\varepsilon b_{0}(ka_{0}+\varepsilon
k(n-1)a_{4}+a_{7}r)R(U,V,W,X)-\varepsilon a_{1}b_{0}S(X,R(U,V)W)  \nonumber
\\
\qquad \qquad \qquad &=&\ -\,2k(n-1)a_{3}(kb_{0}+k(n-1)b_{4}+b_{7}r)\eta
(X)\eta (U)g(V,W)  \nonumber \\
&&-\,2k(n-1)a_{3}(-kb_{0}+k(n-1)b_{5}-b_{7}r)\eta (X)\eta (V)g(U,W) 
\nonumber \\
&&+\,\varepsilon a_{1}b_{4}S^{2}(X,U)g(V,W)+\varepsilon
a_{1}b_{5}S^{2}(X,V)g(U,W)  \nonumber \\
&&+\,\varepsilon a_{1}b_{6}S^{2}(X,W)g(U,V)-a_{5}(b_{1}+b_{3})S^{2}(X,V)\eta
(U)\eta (W)  \nonumber \\
&&-\,a_{5}(b_{1}+b_{2})S^{2}(X,W)\eta (U)\eta
(V)-a_{5}(b_{2}+b_{3})S^{2}(X,U)\eta (V)\eta (W)  \nonumber \\
&&-\,2a_{6}b_{1}S^{2}(V,W)\eta (X)\eta (U)-2a_{6}b_{2}S^{2}(U,W)\eta (X)\eta
(V)  \nonumber \\
&&-\,2a_{6}b_{3}S^{2}(U,V)\eta (X)\eta (W)-2k^{2}(n-1)a_{3}b_{6}g(U,V)\eta
(X)\eta (W)  \nonumber \\
&&-\,2\left( k(n-1)a_{3}b_{1}+a_{6}(kb_{0}+k(n-1)b_{4}+b_{7}r)\right) \eta
(X)\eta (U)S(V,W)  \nonumber \\
&&-\,2\left( k(n-1)a_{3}b_{2}+a_{6}(-kb_{0}+k(n-1)b_{5}-b_{7}r)\right) \eta
(X)\eta (V)S(U,W)  \nonumber \\
&&-\,2k(n-1)(a_{3}b_{3}+a_{6}b_{6})S(U,V)\eta (X)\eta (W)  \nonumber \\
&&+\,\varepsilon
(b_{4}(ka_{0}+k(n-1)a_{4}+a_{7}r)-a_{1}(kb_{0}+k(n-1)b_{4}))S(X,U)g(V,W) 
\nonumber \\
&&+\,\varepsilon
(b_{5}(ka_{0}+k(n-1)a_{4}+a_{7}r)-a_{1}(-kb_{0}+k(n-1)b_{5}))S(X,V)g(U,W) 
\nonumber \\
&&+\,\varepsilon b_{6}(ka_{0}+k(n-1)(a_{4}-a_{1})+a_{7}r)S(X,W)g(U,V) 
\nonumber \\
&&-\,\varepsilon (kb_{0}+k(n-1)b_{4})(ka_{0}+k(n-1)a_{4}+a_{7}r)g(X,U)g(V,W)
\nonumber \\
&&-\,\varepsilon (-kb_{0}+k(n-1)b_{5})(ka_{0}+k(n-1)a_{4}+a_{7}r)g(U,W)g(X,V)
\nonumber \\
&&-\,\varepsilon k(n-1)b_{6}(ka_{0}+k(n-1)a_{4}+a_{7}r)g(X,W)g(U,V) 
\nonumber \\
&&-\,k(n-1)\left( (b_{2}+b_{3})(ka_{0}+k(n-1)a_{4}+a_{7}r)\right.  \nonumber
\\
&&\left. +\,(a_{2}+a_{4})(-kb_{0}+k(n-1)(b_{5}+b_{6})-b_{7}r)\right)
g(X,U)\eta (V)\eta (W)  \nonumber \\
&&-\,k(n-1)\left( (b_{1}+b_{3})(ka_{0}+k(n-1)a_{4}+a_{7}r)\right.  \nonumber
\\
&&\left. +\,(a_{2}+a_{4})(kb_{0}+k(n-1)(b_{4}+b_{6})+b_{7}r)\right)
g(X,V)\eta (U)\eta (W)  \nonumber \\
&&-\,\left( (b_{1}+b_{3})(-ka_{0}+k(n-1)(a_{1}+a_{2})-a_{7}r)\right. 
\nonumber \\
&&\left. +\,(a_{1}+a_{5})(kb_{0}+k(n-1)(b_{4}+b_{6})+b_{7}r)\right)
S(X,V)\eta (U)\eta (W)  \nonumber \\
&&-\,\left( (b_{2}+b_{3})(-ka_{0}+k(n-1)(a_{1}+a_{2})-a_{7}r)\right. 
\nonumber \\
&&\left. +\,(a_{1}+a_{5})(kb_{0}+k(n-1)(b_{5}+b_{6})+b_{7}r)\right)
S(X,U)\eta (V)\eta (W)  \nonumber \\
&&-\,\left( (b_{1}+b_{2})(-ka_{0}+k(n-1)(a_{1}+a_{2})-a_{7}r)\right. 
\nonumber \\
&&\left. +\,k(n-1)(b_{4}+b_{5})(a_{1}+a_{5})\right) S(X,W)\eta (U)\eta (V) 
\nonumber \\
&&-\,k(n-1)\left( k(n-1)(b_{4}+b_{5})(a_{2}+a_{4})\right.  \nonumber \\
&&\left. +\,(b_{1}+b_{2})(ka_{0}+k(n-1)a_{4}+a_{7}r)\right) g(X,W)\eta
(U)\eta (V).  \label{eq-T-T-2}
\end{eqnarray}

\noindent In particular, if $M$ is an $n$-dimensional $({\cal T}_{\!a},{\cal %
T}_{\!a})$-semisymmetric $(N(k),\xi )$-semi-Riemannian manifold, then 
\vspace{-0.8cm} 
\begin{eqnarray*}
&&-\,\varepsilon a_{0}(ka_{0}+\varepsilon
k(n-1)a_{4}+a_{7}r)R(U,V,W,X)-\varepsilon a_{0}a_{1}S(X,R(U,V)W) \\
\qquad \qquad \qquad &=&\ -\,2k(n-1)a_{3}(ka_{0}+k(n-1)a_{4}+a_{7}r)\eta
(X)\eta (U)g(V,W) \\
&&-\,2k(n-1)a_{3}(-ka_{0}+k(n-1)a_{5}-a_{7}r)\eta (X)\eta (V)g(U,W) \\
&&+\,\varepsilon a_{1}a_{4}S^{2}(X,U)g(V,W)+\varepsilon
a_{1}a_{5}S^{2}(X,V)g(U,W) \\
&&+\,\varepsilon a_{1}a_{6}S^{2}(X,W)g(U,V)-a_{5}(a_{1}+a_{3})S^{2}(X,V)\eta
(U)\eta (W) \\
&&-\,a_{5}(a_{1}+a_{2})S^{2}(X,W)\eta (U)\eta
(V)-a_{5}(a_{2}+a_{3})S^{2}(X,U)\eta (V)\eta (W) \\
&&-\,2a_{1}a_{6}S^{2}(V,W)\eta (X)\eta (U)-2a_{2}a_{6}S^{2}(U,W)\eta (X)\eta
(V) \\
&&-\,2a_{3}a_{6}S^{2}(U,V)\eta (X)\eta (W)-2k^{2}(n-1)a_{3}a_{6}g(U,V)\eta
(X)\eta (W) \\
&&-\,2\left( k(n-1)a_{1}a_{3}+a_{6}(ka_{0}+k(n-1)a_{4}+a_{7}r)\right) \eta
(X)\eta (U)S(V,W) \\
&&-\,2\left( k(n-1)a_{2}a_{3}+a_{6}(-ka_{0}+k(n-1)a_{5}-a_{7}r)\right) \eta
(X)\eta (V)S(U,W) \\
&&-\,2k(n-1)(a_{3}a_{3}+a_{6}a_{6})S(U,V)\eta (X)\eta (W) \\
&&+\,\varepsilon
(a_{4}(ka_{0}+k(n-1)a_{4}+a_{7}r)-a_{1}(ka_{0}+k(n-1)a_{4}))S(X,U)g(V,W) \\
&&+\,\varepsilon
(a_{5}(ka_{0}+k(n-1)a_{4}+a_{7}r)-a_{1}(-ka_{0}+k(n-1)a_{5}))S(X,V)g(U,W) \\
&&+\,\varepsilon a_{6}(ka_{0}+k(n-1)(a_{4}-a_{1})+a_{7}r)S(X,W)g(U,V) \\
&&-\,\varepsilon (ka_{0}+k(n-1)a_{4})(ka_{0}+k(n-1)a_{4}+a_{7}r)g(X,U)g(V,W)
\\
&&-\,\varepsilon (-ka_{0}+k(n-1)a_{5})(ka_{0}+k(n-1)a_{4}+a_{7}r)g(U,W)g(X,V)
\\
&&-\,\varepsilon k(n-1)a_{6}(ka_{0}+k(n-1)a_{4}+a_{7}r)g(X,W)g(U,V) \\
&&-\,k(n-1)\left( (a_{2}+a_{3})(ka_{0}+k(n-1)a_{4}+a_{7}r)\right. \\
&&\left. +\,(a_{2}+a_{4})(-ka_{0}+k(n-1)(a_{5}+a_{6})-a_{7}r)\right)
g(X,U)\eta (V)\eta (W) \\
&&-\,k(n-1)\left( (a_{1}+a_{3})(ka_{0}+k(n-1)a_{4}+a_{7}r)\right. \\
&&\left. +\,(a_{2}+a_{4})(ka_{0}+k(n-1)(a_{4}+a_{6})+a_{7}r)\right)
g(X,V)\eta (U)\eta (W) \\
&&-\,\left( (a_{1}+a_{3})(-ka_{0}+k(n-1)(a_{1}+a_{2})-a_{7}r)\right. \\
&&\left. +\,(a_{1}+a_{5})(ka_{0}+k(n-1)(a_{4}+a_{6})+a_{7}r)\right)
S(X,V)\eta (U)\eta (W) \\
&&-\,\left( (a_{2}+a_{3})(-ka_{0}+k(n-1)(a_{1}+a_{2})-a_{7}r)\right. \\
&&\left. +\,(a_{1}+a_{5})(ka_{0}+k(n-1)(a_{5}+a_{6})+a_{7}r)\right)
S(X,U)\eta (V)\eta (W) \\
&&-\,\left( (a_{1}+a_{2})(-ka_{0}+k(n-1)(a_{1}+a_{2})-a_{7}r)\right. \\
&&\left. +\,k(n-1)(a_{4}+a_{5})(a_{1}+a_{5})\right) S(X,W)\eta (U)\eta (V) \\
&&-\,k(n-1)\left( k(n-1)(a_{4}+a_{5})(a_{2}+a_{4})\right. \\
&&\left. +\,(a_{1}+a_{2})(ka_{0}+k(n-1)a_{4}+a_{7}r)\right) g(X,W)\eta
(U)\eta (V).
\end{eqnarray*}
\end{th}

\noindent {\bf Proof.} Let $M$ be an $n$-dimensional $({\cal T}_{\!a},{\cal T%
}_{\!b})$-semisymmetric $(N(k),\xi )$-semi-Riemannian manifold. Then 
\begin{equation}
({\cal T}_{\!a}(Z,X)\cdot {\cal T}_{\!b})(U,V)W=0.  \label{eq-txxi}
\end{equation}%
Taking $Z=\xi $ in (\ref{eq-txxi}), we get 
\[
({\cal T}_{a}(\xi ,X)\cdot {\cal T}_{\!b})(U,V)W=0 
\]%
which gives 
\[
\left[ {\cal T}_{\!a}(\xi ,X),{\cal T}_{\!b}(U,V)\right] W-{\cal T}_{\!b}(%
{\cal T}_{\!a}(\xi ,X)U,V)W-{\cal T}_{\!b}(U,{\cal T}_{\!a}(\xi ,X)V)W=0, 
\]%
that is,%
\begin{eqnarray}
&&0={\cal T}_{\!a}(\xi ,X){\cal T}_{\!b}(U,V)W-{\cal T}_{\!b}({\cal T}%
_{\!a}(\xi ,X)U,V)W  \nonumber \\
&&\quad -\ {\cal T}_{\!b}(U,{\cal T}_{\!a}(\xi ,X)V)W-{\cal T}_{\!b}(U,V)%
{\cal T}_{\!a}(\xi ,X)W.  \label{eq-T-T}
\end{eqnarray}%
Taking the inner product of (\ref{eq-T-T}) with $\xi $, we get 
\begin{eqnarray}
&&0={\cal T}_{\!a}(\xi ,X,{\cal T}_{\!b}(U,V)W,\xi )-{\cal T}_{\!b}({\cal T}%
_{\!a}(\xi ,X)U,V,W,\xi )  \nonumber \\
&&\quad -\ {\cal T}_{\!b}(U,{\cal T}_{\!a}(\xi ,X)V,W,\xi )-{\cal T}%
_{\!b}(U,V,{\cal T}_{\!a}(\xi ,X)W,\xi ).  \label{eq-T-T-1}
\end{eqnarray}%
By using (\ref{eq-X-Y-xi}),\ldots ,(\ref{eq-X-xi-xi}) in (\ref{eq-T-T-1}),
we get (\ref{eq-T-T-2}). $\blacksquare $

\begin{th}
Let $M$ be an $n$-dimensional $({\cal T}_{\!a},{\cal T}_{\!b})$%
-semisymmetric $(N(k),\xi )$-semi-Riemannian manifold. Then 
\begin{eqnarray}
&&\varepsilon \left\{
a_{5}b_{0}+na_{5}b_{1}+a_{5}b_{2}+a_{5}b_{6}+a_{5}b_{3}+a_{5}b_{5}\right\}
S^{2}(V,W)  \nonumber \\
&=&\ (\left\{ (nb_{1}+b_{2}+b_{3}+b_{5}+b_{6}+b_{0})(\varepsilon
ka_{0}+\varepsilon b_{7}r)\right.  \nonumber \\
&&\qquad -\ \varepsilon
k(n-1)(2a_{5}b_{6}+a_{2}b_{3}+a_{1}b_{6}+a_{1}b_{3}+a_{1}b_{5}  \nonumber \\
&&\qquad \left. +\
a_{1}b_{1}+a_{1}b_{2}+a_{2}b_{2}+a_{2}b_{6}+na_{2}b_{1}+a_{1}b_{0}+a_{2}b_{0})\right\}
\nonumber \\
&&-\ \varepsilon (n-1)a_{1}b_{7}r-\varepsilon na_{5}b_{7}r-\varepsilon
b_{4}a_{5}r-\varepsilon a_{1}b_{4}r)S(V,W)  \nonumber \\
&&+\ \left\{ -\varepsilon
k(n-1)(nb_{1}+b_{2}+b_{3}+b_{5}+b_{6}+b_{0})(a_{7}r+ka_{0}+k(n-1)a_{4})%
\right.  \nonumber \\
&&\qquad \left. -\ \varepsilon
k(n-1)r((n-1)b_{7}a_{2}+(n-1)b_{7}a_{4}+a_{2}b_{4}+a_{4}b_{4})\right\} g(V,W)
\nonumber \\
&&+\ (a_{1}+a_{2}+2a_{3}+a_{4}+a_{5}+2a_{6})\left\{
-k^{2}(n-1)^{2}(nb_{1}+b_{2}+b_{3}+b_{5}+b_{6}+b_{0})\right.  \nonumber \\
&&\qquad \left. -\ k(n-1)^{2}b_{7}r-k(n-1)b_{4}r\right\} \eta (V)\eta (W).
\label{eq-semi-TS}
\end{eqnarray}%
In particular, if $M$ is an $n$-dimensional $({\cal T}_{\!a},{\cal T}_{\!a})$%
-semisymmetric $(N(k),\xi )$-semi-Riemannian manifold, then 
\begin{eqnarray*}
&&\varepsilon \left\{
a_{5}a_{0}+na_{5}a_{1}+a_{5}a_{2}+a_{5}a_{6}+a_{5}a_{3}+a_{5}^{2}\right\}
S^{2}(V,W) \\
&=&\ (\left\{ (na_{1}+a_{2}+a_{3}+a_{5}+a_{6}+a_{0})(\varepsilon
ka_{0}+\varepsilon a_{7}r)\right. \\
&&\qquad -\ \varepsilon
k(n-1)(2a_{5}a_{6}+a_{2}a_{3}+a_{1}a_{6}+a_{1}a_{3}+a_{1}a_{5} \\
&&\qquad \left. +\
a_{1}^{2}+a_{1}a_{2}+a_{2}^{2}+a_{2}a_{6}+na_{1}a_{2}+a_{0}a_{1}+a_{0}a_{2})%
\right\} \\
&&-\ \varepsilon (n-1)a_{1}a_{7}r-\varepsilon na_{5}a_{7}r-\varepsilon
a_{4}a_{5}r-\varepsilon a_{1}a_{4}r)S(V,W) \\
&&+\ \left\{ -\varepsilon
k(n-1)(na_{1}+a_{2}+a_{3}+a_{5}+a_{6}+a_{0})(a_{7}r+ka_{0}+k(n-1)a_{4})%
\right. \\
&&\qquad \left. -\ \varepsilon
k(n-1)r((n-1)a_{2}a_{7}+(n-1)a_{4}a_{7}+a_{2}a_{4}+a_{4}^{2})\right\} g(V,W)
\\
&&+\ (a_{1}+a_{2}+2a_{3}+a_{4}+a_{5}+2a_{6})\left\{
-k^{2}(n-1)^{2}(na_{1}+a_{2}+a_{3}+a_{5}+a_{6}+a_{0})\right. \\
&&\qquad \left. -\ k(n-1)^{2}a_{7}r-k(n-1)a_{4}r\right\} \eta (V)\eta (W).
\end{eqnarray*}
\end{th}

\noindent {\bf Proof.} By contracting (\ref{eq-T-T-2}) we get (\ref%
{eq-semi-TS}). $\blacksquare $ \medskip

From Theorem \ref{th-T-T}, we get

\begin{th}
\label{GCT-ss} Let $M$ be an $n$-dimensional ${\cal T}_{\!a}$-semisymmetric $%
(N(k),\xi )$-semi-Riemannian manifold. Then 
\begin{eqnarray}
-\varepsilon a_{0}kR(U,V,W,X) &=&\varepsilon ka_{4}S(X,U)g(V,W)+\varepsilon
ka_{5}S(X,V)g(U,W)  \nonumber \\
&&+\,\varepsilon ka_{6}S(X,W)g(U,V)-\varepsilon k^{2}(n-1)a_{6}g(X,W)g(U,V) 
\nonumber \\
&&-\,\varepsilon k(ka_{0}+k(n-1)a_{4})g(V,W)g(X,U)  \nonumber \\
&&-\,\varepsilon k(-ka_{0}+k(n-1)a_{5})g(U,W)g(X,V)  \nonumber \\
&&-\,k^{2}(n-1)(a_{2}+a_{3})g(X,U)\eta (V)\eta (W)  \nonumber \\
&&-\,k^{2}(n-1)(a_{1}+a_{3})g(X,V)\eta (U)\eta (W)  \nonumber \\
&&-\,k^{2}(n-1)(a_{1}+a_{2})g(X,W)\eta (U)\eta (V)  \nonumber \\
&&+\,k(a_{2}+a_{3})S(X,U)\eta (V)\eta (W)  \nonumber \\
&&+\,k(a_{1}+a_{3})S(X,V)\eta (U)\eta (W)  \nonumber \\
&&+\,k(a_{1}+a_{2})S(X,W)\eta (U)\eta (V).  \label{eq-semi-sym-R}
\end{eqnarray}
\end{th}

\begin{th}
\label{GCT-sss} Let $M$ be an $n$-dimensional ${\cal T}_{\!a}$-semisymmetric 
$(N(k),\xi )$-semi-Riemannian manifold such that $a_{0}+a_{5}+a_{6}\not=0$.
Then

\begin{description}
\item[{\bf (a)}] 
\begin{equation}
S(V,W)=B_{1}g(V,W)+B_{2}\eta (V)\eta (W)  \label{eq-SVW}
\end{equation}
and 
\begin{eqnarray}
-a_{0}R(U,V,W,X) &=&\left( a_{4}B_{1}-ka_{0}-k(n-1)a_{4}\right) g(X,U)g(V,W)
\nonumber \\
&&+\,\left( a_{5}B_{1}+ka_{0}-k(n-1)a_{5}\right) g(X,V)g(U,W)  \nonumber \\
&&+\,\left( a_{6}B_{1}-k(n-1)a_{6}\right) g(X,W)g(U,V)  \nonumber \\
&&+\,\varepsilon (a_{2}+a_{3})\left( B_{1}-k(n-1)\right) g(X,U)\eta (V)\eta
(W)  \nonumber \\
&&+\,\varepsilon (a_{1}+a_{3})\left( B_{1}-k(n-1)\right) g(X,V)\eta (U)\eta
(W)  \nonumber \\
&&+\,\varepsilon (a_{1}+a_{2})\left( B_{1}-k(n-1)\right) g(X,W)\eta (U)\eta
(V)  \nonumber \\
&&+\,2\varepsilon B_{2}(a_{1}+a_{2}+a_{3})\eta (X)\eta (U)\eta (V)\eta (W) 
\nonumber \\
&&+\,a_{4}B_{2}g(V,W)\eta (X)\eta (U)  \nonumber \\
&&+\,a_{5}B_{2}g(U,W)\eta (X)\eta (V)  \nonumber \\
&&+\,a_{6}B_{2}g(U,V)\eta (X)\eta (W),  \label{eq-RSVW}
\end{eqnarray}
where 
\[
B_{1}=-\,\frac{a_{4}r-n(ka_{0}+k(n-1)a_{4})-(-ka_{0}+k(n-1)a_{5})-k(n-1)a_{6}%
}{a_{0}+a_{5}+a_{6}}, 
\]
\[
B_{2}=-\,\frac{\varepsilon k(a_{2}+a_{3})(r-n(n-1)k)}{a_{0}+a_{5}+a_{6}}. 
\]

\item[{\bf (b)}] If $a_{0}+a_{2}+a_{3}+na_{4}+a_{5}+a_{6}\not=0$, then it is
an Einstein manifold and a manifold of constant curvature $k$.
\end{description}
\end{th}

\noindent {\bf Proof. }Let $M$ be an $n$-dimensional ${\cal T}_{\!a}$%
-semisymmetric $(N(k),\xi )$-semi-Riemannian manifold such that $%
a_{0}+a_{5}+a_{6}\not=0$.\newline
\noindent {\bf Case (a).} Contracting (\ref{eq-semi-sym-R}), we get (\ref%
{eq-SVW}). \newline
\noindent {\bf Case (b).} Let $a_{0}+a_{2}+a_{3}+na_{4}+a_{5}+a_{6}\not=0$,
then contracting (\ref{eq-SVW}), we get 
\begin{equation}
r=\,kn(n-1).  \label{eq-gen-r}
\end{equation}%
Since $a_{0}+a_{5}+a_{6}\not=0$, then by (\ref{eq-semi-sym-R}) and (\ref%
{eq-gen-r}) we get 
\begin{equation}
S(V,W)=\,k(n-1)g(V,W).  \label{eq-gen-S}
\end{equation}%
Using (\ref{eq-gen-S}) and (\ref{eq-gen-r}) in (\ref{eq-semi-sym-R}), we get%
\begin{equation}
R(U,V,W,X)=k\left( g(V,W)g(X,U)-g(U,W)g(X,V\right) ),  \label{eq-R-111}
\end{equation}%
which proves the result. $\blacksquare $

In view of Theorem~\ref{GCT-sss}, we have the following

\begin{cor}
\label{cor-1} Let $M$ be an $n$-dimensional ${\cal T}_{\!a}$-semisymmetric $%
\left( N(k),\xi \right) $-semi-Riemannian manifold such that 
\[
{\cal T}_{\!a}\in \left\{ R,{\cal V},{\cal P},{\cal M},{\cal W}_{0},{\cal W}%
_{0}^{\ast },{\cal W}_{1},{\cal W}_{1}^{\ast },{\cal W}_{3},\ldots ,{\cal W}%
_{8}\right\} 
\]%
Then we have the following two tables\/{\rm :}~%
\[
\begin{tabular}{|l|l|}
\hline
${\boldmath M}$ & ${\boldmath S=}$ \\ \hline
$N(k)$-contact metric & $k(n-1)g$ \\ \hline
Sasakian & $(n-1)g$ \\ \hline
Kenmotsu & $-\,(n-1)g$ \\ \hline
$(\varepsilon )$-Sasakian & $\varepsilon (n-1)g$ \\ \hline
para-Sasakian & $-\,(n-1)g$ \\ \hline
$(\varepsilon )$-para-Sasakian & $-\,\varepsilon (n-1)g$ \\ \hline
\end{tabular}%
\]%
\[
\begin{tabular}{|l|l|}
\hline
${\boldmath M}$ & ${\boldmath R(X,Y)Z=}$ \\ \hline
$N(k)$-contact metric & $k\left( g(Y,Z)X-g(X,Z)Y\right) $ \\ \hline
Sasakian & $g(Y,Z)X-g(X,Z)Y$ \\ \hline
Kenmotsu & $-\,\left( g(Y,Z)X-g(X,Z)Y\right) $ \\ \hline
$(\varepsilon )$-Sasakian & $\varepsilon \left( g(Y,Z)X-g(X,Z)Y\right) $ \\ 
\hline
para-Sasakian & $-\,\left( g(Y,Z)X-g(X,Z)Y\right) $ \\ \hline
$(\varepsilon )$-para-Sasakian & $-\,\varepsilon \left(
g(Y,Z)X-g(X,Z)Y\right) $ \\ \hline
\end{tabular}%
\]
\end{cor}

\begin{cor}
\label{cor-2} Let $M$ be an $n$-dimensional quasi-conformal semisymmetric $%
\left( N(k),\xi \right) $-semi-Riemannian manifold such that $%
a_{0}-a_{1}\not=0$\ and $a_{0}+(n-2)a_{1}\not=0$. Then we have the following
two tables\/{\rm :}~%
\[
\begin{tabular}{|l|l|}
\hline
${\boldmath M}$ & ${\boldmath S=}$ \\ \hline
$N(k)$-contact metric & $k(n-1)g$ \\ \hline
Sasakian {\rm \cite{De-Jun-Gazi-08}} & $(n-1)g$ \\ \hline
Kenmotsu & $-\,(n-1)g$ \\ \hline
$(\varepsilon )$-Sasakian & $\varepsilon (n-1)g$ \\ \hline
para-Sasakian & $-\,(n-1)g$ \\ \hline
$(\varepsilon )$-para-Sasakian & $-\,\varepsilon (n-1)g$ \\ \hline
\end{tabular}%
\]%
\[
\begin{tabular}{|l|l|}
\hline
${\boldmath M}$ & ${\boldmath R(X,Y)Z=}$ \\ \hline
$N(k)$-contact metric & $k\left( g(Y,Z)X-g(X,Z)Y\right) $ \\ \hline
Sasakian & $g(Y,Z)X-g(X,Z)Y$ \\ \hline
Kenmotsu & $-\,\left( g(Y,Z)X-g(X,Z)Y\right) $ \\ \hline
$(\varepsilon )$-Sasakian & $\varepsilon \left( g(Y,Z)X-g(X,Z)Y\right) $ \\ 
\hline
para-Sasakian & $-\,\left( g(Y,Z)X-g(X,Z)Y\right) $ \\ \hline
$(\varepsilon )$-para-Sasakian & $-\,\varepsilon \left(
g(Y,Z)X-g(X,Z)Y\right) $ \\ \hline
\end{tabular}%
\]
\end{cor}

\begin{cor}
\label{cor-3} Let $M$ be an $n$-dimensional pseudo-projective semisymmetric $%
\left( N(k),\xi \right) $-semi-Riemannian manifold such that $a_{0}\not=0$
and $a_{0}-a_{1}\not=0$. Then we have the following two tables\/{\rm :}~ 
\[
\begin{tabular}{|l|l|}
\hline
${\boldmath M}$ & ${\boldmath S=}$ \\ \hline
$N(k)$-contact metric & $k(n-1)g$ \\ \hline
Sasakian & $(n-1)g$ \\ \hline
Kenmotsu & $-\,(n-1)g$ \\ \hline
$(\varepsilon )$-Sasakian & $\varepsilon (n-1)g$ \\ \hline
para-Sasakian & $-\,(n-1)g$ \\ \hline
$(\varepsilon )$-para-Sasakian & $-\,\varepsilon (n-1)g$ \\ \hline
\end{tabular}%
\]%
\[
\begin{tabular}{|l|l|}
\hline
${\boldmath M}$ & ${\boldmath R(X,Y)Z=}$ \\ \hline
$N(k)$-contact metric & $k\left( g(Y,Z)X-g(X,Z)Y\right) $ \\ \hline
Sasakian & $g(Y,Z)X-g(X,Z)Y$ \\ \hline
Kenmotsu & $-\,\left( g(Y,Z)X-g(X,Z)Y\right) $ \\ \hline
$(\varepsilon )$-Sasakian & $\varepsilon \left( g(Y,Z)X-g(X,Z)Y\right) $ \\ 
\hline
para-Sasakian & $-\,\left( g(Y,Z)X-g(X,Z)Y\right) $ \\ \hline
$(\varepsilon )$-para-Sasakian & $-\,\varepsilon \left(
g(Y,Z)X-g(X,Z)Y\right) $ \\ \hline
\end{tabular}%
\]
\end{cor}

Here, we give the well known results of Okumura \cite{Okumura-62} and
Koufogiorgos \cite{Koufogiorgos-93}.

\begin{th}
{\rm \cite[Lemma 2.2]{Okumura-62} }If an $n$-dimensional Sasakian manifold
is conformally flat, then the scalar curvature has a positive constant value 
$n(n-1)$.
\end{th}

\begin{th}
{\rm \cite[Corollary 3.3]{Koufogiorgos-93}} Let $M$ be an $\eta $-Einstein
contact metric manifold of dimension $2m+1>5$. If $\xi $ belongs to the $k$%
-nullity distribution, then $k=1$ and the structure is Sasakian.
\end{th}

\begin{cor}
\label{cor-4} Let $M$ be an $n$-dimensional Weyl-semisymmetric $\left(
N(k),\xi \right) $-semi-Riemannian manifold. Then we have the following two
tables\/{\rm :} 
\[
\begin{tabular}{|l|l|}
\hline
${\boldmath M}$ & ${\boldmath S=}$ \\ \hline
$N(k)$-contact metric & $k(n-1)g$ \\ \hline
Sasakian & $(n-1)g$ \\ \hline
Kenmotsu & $-\,(n-1)g$ \\ \hline
$(\varepsilon )$-Sasakian & $\varepsilon (n-1)g$ \\ \hline
para-Sasakian {\rm \cite{Tarafdar-Mayra-94}} & $-\,(n-1)g$ \\ \hline
$(\varepsilon )$-para-Sasakian & $-\,\varepsilon (n-1)g$ \\ \hline
\end{tabular}%
\]%
\[
\begin{tabular}{|l|l|}
\hline
${\boldmath M}$ & ${\boldmath R(X,Y)Z=}$ \\ \hline
$N(k)$-contact metric & $k\left( g(Y,Z)X-g(X,Z)Y\right) $ \\ \hline
Sasakian {\rm \cite{Chaki-Tarafdar-90}} & $g(Y,Z)X-g(X,Z)Y$ \\ \hline
Kenmotsu & $-\,\left( g(Y,Z)X-g(X,Z)Y\right) $ \\ \hline
$(\varepsilon )$-Sasakian & $\varepsilon \left( g(Y,Z)X-g(X,Z)Y\right) $ \\ 
\hline
para-Sasakian & $-\,\left( g(Y,Z)X-g(X,Z)Y\right) $ \\ \hline
$(\varepsilon )$-para-Sasakian & $-\,\varepsilon \left(
g(Y,Z)X-g(X,Z)Y\right) $ \\ \hline
\end{tabular}%
\]
\end{cor}

\noindent {\bf Proof.} By putting the value for conformal curvature tensor
in (\ref{eq-SVW}), we get 
\begin{equation}
S=\left( \frac{r}{n-1}-k\right) g+\varepsilon k\left( nk-\frac{r}{n-1}%
\right) \eta \otimes \eta .  \label{eq-SVWW}
\end{equation}%
{\bf Case 1.} Let $k\not=1$. Contracting (\ref{eq-SVWW}), we get 
\[
r=kn(n-1),\quad k\not=1. 
\]%
Using the value of $r$ in (\ref{eq-SVWW}), we get 
\[
S=k(n-1)g. 
\]%
Using this in (\ref{eq-RSVW}), we get 
\[
R(X,Y)Z=k\left( g(Y,Z)X-g(X,Z)Y\right) . 
\]%
{\bf Case 2. }Let $k=1$. By putting the value for conformal curvature tensor
in (\ref{eq-RSVW}), we get ${\cal C}=0$. Then using the result of \cite%
{Okumura-62} and \cite{Koufogiorgos-93}, we get $r=n(n-1)$. $\blacksquare $

\begin{cor}
\label{cor-5} Let $M$ be an $n$-dimensional conharmonic semisymmetric $%
\left( N(k),\xi \right) $-semi-Riemannian manifold. Then we have the
following two tables\/{\rm :}~ 
\[
\begin{tabular}{|l|l|}
\hline
${\boldmath M}$ & ${\boldmath S=}$ \\ \hline
$N(k)$-contact metric & $k(n-1)g$ \\ \hline
Sasakian & $\left( \dfrac{r}{n-1}-1\right) g+\left( n-\dfrac{r}{n-1}\right)
\eta \otimes \eta $ \\ \hline
Kenmotsu & $-\,(n-1)g$ \\ \hline
$(\varepsilon )$-Sasakian & $\varepsilon (n-1)g$ \\ \hline
para-Sasakian & $-\,(n-1)g$ \\ \hline
$(\varepsilon )$-para-Sasakian & $-\,\varepsilon (n-1)g$ \\ \hline
\end{tabular}%
\]%
\[
\begin{tabular}{|l|l|}
\hline
${\boldmath M}$ & ${\boldmath R(X,Y)Z=}$ \\ \hline
$N(k)$-contact metric & $k\left( g(Y,Z)X-g(X,Z)Y\right) $ \\ \hline
Sasakian & $\dfrac{1}{n-2}\left( \dfrac{r}{n-1}-2\right) \left(
g(Y,Z)X-g(X,Z)Y\right) $ \\ 
& $+\,\dfrac{1}{n-2}\left( n-\dfrac{r}{n-1}\right) (\eta (Y)\eta (Z)X-\eta
(X)\eta (Z)Y$ \\ 
& $+\,g(Y,Z)\eta (X)\xi -g(X,Z)\eta (Y)\xi )$ \\ \hline
Kenmotsu & $-\,\left( g(Y,Z)X-g(X,Z)Y\right) $ \\ \hline
$(\varepsilon )$-Sasakian & $\varepsilon \left( g(Y,Z)X-g(X,Z)Y\right) $ \\ 
\hline
para-Sasakian & $-\,\left( g(Y,Z)X-g(X,Z)Y\right) $ \\ \hline
$(\varepsilon )$-para-Sasakian & $-\,\varepsilon \left(
g(Y,Z)X-g(X,Z)Y\right) $ \\ \hline
\end{tabular}%
\]
\end{cor}

\noindent {\bf Proof.} By putting the value for conharmonic curvature tensor
in (\ref{eq-SVW}), we get 
\begin{equation}
S=\left( \frac{r}{n-1}-k\right)g + \varepsilon k\left( nk-\frac{r}{n-1}%
\right) \eta \otimes \eta .  \label{eq-SVWW-1}
\end{equation}%
Let $k\not=1$. Contracting (\ref{eq-SVWW-1}), we get 
\[
r=kn(n-1),\quad k\not=1. 
\]%
Using the value of $r$ in (\ref{eq-SVWW-1}), we get 
\[
S=k(n-1)g. 
\]%
Using this in (\ref{eq-RSVW}), we get 
\[
R(X,Y)Z=k\left( g(Y,Z)X-g(X,Z)Y\right) . 
\]

\begin{cor}
\label{cor-7} Let $M$ be an $n$-dimensional ${\cal W}_{2}$-semisymmetric $%
\left( N(k),\xi \right) $-semi-Riemannian manifold. If $M$ is one of $N(k)$%
-contact metric manifold, Sasakian manifold, Kenmotsu manifold, $%
(\varepsilon )$-Sasakian manifold, para-Sasakian manifold or $(\varepsilon ) 
$-para-Sasakian manifold, then 
\[
S=\frac{r}{n}g 
\]%
and 
\[
R(X,Y)Z=\dfrac{r}{n(n-1)}\left( g(Y,Z)X-g(X,Z)Y\right) . 
\]
\end{cor}

\noindent {\bf Proof.} By putting the values for ${\cal W}_{2}$-curvature
tensor in (\ref{eq-SVW}) and (\ref{eq-RSVW}), we get the result. $%
\blacksquare $

\begin{cor}
\label{cor-6} Let $M$ be an $n$-dimensional ${\cal W}_{9}$-semisymmetric $%
\left( N(k),\xi \right) $-semi-Riemannian manifold. Then we have the
following two tables\/{\rm :}~ 
\[
\begin{tabular}{|l|l|}
\hline
${\boldmath M}$ & ${\boldmath S=}$ \\ \hline
$N(k)$-contact metric & $k(n-1)g$ \\ \hline
Sasakian & $\left( \dfrac{r}{n-1}-1\right) g+\left( n-\dfrac{r}{n-1}\right)
\eta \otimes \eta $ \\ \hline
Kenmotsu & $-\,(n-1)g$ \\ \hline
$(\varepsilon )$-Sasakian & $\varepsilon (n-1)g$ \\ \hline
para-Sasakian & $-\,(n-1)g$ \\ \hline
$(\varepsilon )$-para-Sasakian & $-\,\varepsilon (n-1)g$ \\ \hline
\end{tabular}%
\]%
\[
\begin{tabular}{|l|l|}
\hline
${\boldmath M}$ & ${\boldmath R(X,Y)Z=}$ \\ \hline
$N(k)$-contact metric & $k\left( g(Y,Z)X-g(X,Z)Y\right) $ \\ \hline
Sasakian & $\dfrac{1}{n-1}\left( \dfrac{r}{n-1}-1\right) (g(Y,Z)X-g(X,Z)Y)$
\\ 
& $+\,\dfrac{1}{n-1}\left( n-\dfrac{r}{n-1}\right) (\eta (Y)\eta (Z)X+\eta
(X)\eta (Z)Y)$ \\ 
& $-\,\dfrac{2}{n-1}\left( n-\dfrac{r}{n-1}\right) \eta (X)\eta (Y)\eta
(Z)\xi $ \\ 
& $+\,\dfrac{1}{n-1}\left( n-\dfrac{r}{n-1}\right) g(Y,Z)\eta (X)\xi $ \\ 
\hline
Kenmotsu & $-\,\left( g(Y,Z)X-g(X,Z)Y\right) $ \\ \hline
$(\varepsilon )$-Sasakian & $\varepsilon \left( g(Y,Z)X-g(X,Z)Y\right) $ \\ 
\hline
para-Sasakian & $-\,\left( g(Y,Z)X-g(X,Z)Y\right) $ \\ \hline
$(\varepsilon )$-para-Sasakian & $-\,\varepsilon \left(
g(Y,Z)X-g(X,Z)Y\right) $ \\ \hline
\end{tabular}%
\]
\end{cor}

\noindent {\bf Proof.} By putting the value for ${\cal W}_{9}$-curvature
tensor in (\ref{eq-SVW}), we get 
\begin{equation}
S=\left( \frac{r}{n-1}-k\right) g+\varepsilon k\left( nk-\frac{r}{n-1}%
\right) \eta \otimes \eta .  \label{eq-SVWW-2}
\end{equation}%
Let $k\not=1$. Contracting (\ref{eq-SVWW-2}), we get 
\[
r=kn(n-1),\quad k\not=1. 
\]%
Using the value of $r$ in (\ref{eq-SVWW-2}), we get 
\[
S=k(n-1)g. 
\]%
Using this in (\ref{eq-RSVW}), we get 
\[
R(X,Y)Z=k\left( g(Y,Z)X-g(X,Z)Y\right) . 
\]

\begin{rem-new}
By Theorem {\rm \ref{GCT-re}}, we conclude that the same results hold if the
condition of ${\cal T}_{\!a}$-semisymmetric is replaced by ${\cal T}_{\!a}$%
-recurrent or ${\cal T}_{\!a}$-symmetric.
\end{rem-new}

\begin{rem-new}
Some of the results related to the above Corollaries have been proved by
authors Miyazawa and Yamaguchi (\cite{Miyazawa-66},\cite%
{Miyazawa-Yamaguchi-66}), Mishra \cite{Mishra-70}, Adati and Matsumoto \cite%
{Adati-Matsumoto-77}, Sato and Matsumoto \cite{Sato-Matsumoto-79}, Adati and
Miyazawa \cite{Adati-Miyazawa-79}, Maralabhavi \cite{Maralabhavi-86}, Ojha 
\cite{Ojha-86}, De and Ghosh \cite{De-Ghosh-94}, Rahman \cite{Rahman-95},
Ghosh and Sharma \cite{Ghosh-Sharma-97}, Mishra and Ojha \cite%
{Mishra-Ojha-01}, Tarafdar and Sengupta \cite{Tarafdar-Sengupta-01}, Blair
et al. \cite{Blair-Kim-Tripathi-05}, Jun et al. \cite{Jun-De-Pathak-05}, 
\"{O}zg\"{u}r and De \cite{Ozgur-De-06}, Tripathi et al. \cite{TKYK-09}.
\end{rem-new}

\begin{rem-new}
There are 400 combinations of derivations for the 20 curvature tensors
mentioned as particular cases in Definition~\ref{defn-GCT}. Here, we
discussed the results for 6 different structures for each derivation
condition. So there are total 2400 results for different structures and
curvature tensors. Out of these 2400 cases, we have discussed only 120 cases
in this paper. The remaining 2280 cases can be obtained by putting the
appropriate value for the curvature tensors, $\varepsilon $ and $k$ in (\ref%
{eq-T-T-2}) and (\ref{eq-semi-TS}). Out of the remaining 2280 cases, some
are mentioned below.
\end{rem-new}

\begin{cor}
{\rm \cite{Hong-Ozgur-Tripathi-06}} A $(2m+1)$-dimensional Kenmotsu manifold 
$M$ satisfies ${\cal V}(\xi ,X)\cdot R=0$ if and only if $M$ is either of
constant scalar curvature or of constant curvature $-1$.
\end{cor}

\begin{cor}
{\rm \cite{Hong-Ozgur-Tripathi-06}} A $(2m+1)$-dimensional Kenmotsu manifold 
$M$ satisfies ${\cal V}(\xi ,X)\cdot {\cal V}=0$ if and only if $M$ is
either of constant scalar curvature or of constant curvature $-1$.
\end{cor}

\begin{cor}
{\rm \cite{Ozgur-Tripathi-2007}} An $n$-dimensional para-Sasakian manifold $%
M $ satisfies ${\cal V}(\xi ,X)\cdot {\cal V}=0$ if and only if either the
scalar curvature $r$ of $M$ is $r=n(1-n)$ or $M$ is locally isometric to the
Hyperbolic space $H^{n}(-1)$.
\end{cor}

\begin{cor}
{\rm \cite{Ozgur-Tripathi-2007}} An $n$-dimensional para-Sasakian manifold $%
M $ satisfies ${\cal V}(\xi ,X)\cdot R=0$ if and only if either $M$ is
locally isometric to the Hyperbolic space $H^{n}(-1)$ or $M$ has constant
scalar curvature $r=n(1-n)$.
\end{cor}

\begin{cor}
{\rm \cite{Ozgur-Tripathi-2007}} An $n$-dimensional para-Sasakian manifold $%
M $ satisfies ${\cal V}(\xi ,X)\cdot {\cal C}=0$ if and only if either $M$
has scalar curvature $r=n(1-n)$ or $M$ is conformally flat, in which case $M$
is a special para-Sasakian manifold.
\end{cor}

\begin{cor}
{\rm \cite{Blair-Kim-Tripathi-05}} A $(2m+1)$-dimensional $N(k)$-contact
metric manifold $M$ satisfies ${\cal V}(\xi ,X)\cdot {\cal V}=0$ if and only
if $M$ is locally isometric to the sphere $S^{2m+1}(1)$ or $M$ is $3 $%
-dimensional and flat.
\end{cor}

\begin{cor}
{\rm \cite{Blair-Kim-Tripathi-05}} A $(2m+1)$-dimensional $N(k)$-contact
metric manifold $M$ satisfies ${\cal V}(\xi ,X)\cdot R=0$ if and only if $M$
is locally isometric to the sphere $S^{2m+1}(1)$ or $M$ is $3$-dimensional
and flat.
\end{cor}

\begin{th}
Let $M$ be an $n$-dimensional ${\cal T}_{\!a}$-semisymmetric $\left(
N(k),\xi \right) $-semi-Riemannian manifold such that $a_{0}+a_{5}+a_{6}%
\not=0$ and 
\[
a_{0}+a_{2}+a_{3}+na_{4}+a_{5}+a_{6}\not=0, 
\]%
then $M$ is ${\cal T}_{\!a}$-flat if either $k=0$ or $k\not=0$ and 
\[
a_{0}+a_{1}(n-1)+a_{4}(n-1)+a_{7}n(n-1)=0, 
\]%
\[
-a_{0}+a_{2}(n-1)+a_{5}(n-1)-a_{7}n(n-1)=0, 
\]%
\[
a_{3}+a_{6}=0. 
\]
\end{th}

\begin{th}
Let $M$ be an $n$-dimensional $\left( N(k),\xi \right) $-semi-Riemannian
manifold of constant curvature $k$ and 
\[
a_{0}+a_{1}(n-1)+a_{4}(n-1)+a_{7}n(n-1)=0, 
\]%
\[
-a_{0}+a_{2}(n-1)+a_{5}(n-1)-a_{7}n(n-1)=0, 
\]%
\[
a_{3}+a_{6}=0, 
\]%
then it is ${\cal T}_{\!a}$-semisymmetric.
\end{th}

\begin{cor}
Let $M$ be an $\left( N(k),\xi \right) $-semi-Riemannian manifold such that 
\[
{\cal T}_{\!a}\in \left\{ {\cal C}_{\ast },{\cal C},{\cal V},{\cal P}_{\ast
},{\cal P},{\cal M},{\cal W}_{0},{\cal W}_{1}^{\ast },{\cal W}_{2}\right\} . 
\]%
Then it is ${\cal T}_{\!a}$-semisymmetric if and only if it is a manifold of
constant curvature $k$.
\end{cor}

\begin{th}
Let $M$ be ${\cal T}_{\!a}$-semisymmetric $\left( N(k),\xi \right) $%
-semi-Riemannian manifold such that $a_{0}+a_{5}+a_{6}\not=0$ and 
\[
a_{0}+a_{2}+a_{3}+na_{4}+a_{5}+a_{6}\not=0, 
\]%
then $M$ is ${\cal T}_{\!a}$-conservative.
\end{th}

\noindent {\bf Proof. }Using (\ref{eq-gen-S}) and (\ref{eq-gen-r}) in (\ref%
{eq-divT-11}), we get {\rm div }${\cal T}_{\!a}=0$.\ $\blacksquare $

\begin{ex}
{\rm \cite{Koufogiorgos-97}} Let $M$ be an $n$-dimensional $N(k)$-contact
metric manifold. If the $\varphi $-sectional curvature of any point of $M$
is independent of the choice of $\varphi $-section at the point, then it is
constant on $M$ and the curvature tensor is given by 
\begin{eqnarray}
4R(X,Y)Z &=&(c+3)(g(Y,Z)X-g(X,Z)Y)  \nonumber \\
&&+(c+3-4k)(\eta (X)\eta (Z)Y-\eta (Y)\eta (Z)X  \nonumber \\
&&+g(X,Z)\eta (Y)\xi -g(Y,Z)\eta (X)\xi )  \nonumber \\
&&+(c-1)(2g(X,\varphi Y)\varphi Z+g(X,\varphi Z)\varphi Y  \nonumber \\
&&-g(Y,\varphi Z)\varphi X),  \label{eq-H-1}
\end{eqnarray}%
where $c$ is the constant $\varphi $-sectional curvature.

Contracting $(\ref{eq-H-1})$, we get 
\begin{equation}
4S(Y,Z)=E_{1}g(Y,Z)-E_{2}\eta (Y)\eta (Z),  \label{eq-H-2}
\end{equation}%
where 
\[
E_{1}=(n-1)(c+3)-(c+3-4k)+3(c-1) 
\]%
and 
\[
E_{2}=(n-2)(c+3-4k)+3(c-1). 
\]%
Consider a ${\cal T}_{\!a}$-semisymmetric $N(k)$-contact metric manifold $M$
with constant $\varphi $-sectional curvature $c$ such that 
\[
{\cal T}_{\!a}\in \left\{ R,{\cal C}_{\ast },{\cal C},{\cal L},{\cal V},%
{\cal P}_{\ast },{\cal P},{\cal M},{\cal W}_{0},{\cal W}_{0}^{\ast },{\cal W}%
_{1},{\cal W}_{1}^{\ast },{\cal W}_{3},\ldots ,{\cal W}_{9}\right\} , 
\]%
then by Corollaries $\ref{cor-1}$, $\ref{cor-2}$, $\ref{cor-3}$, $\ref{cor-4}
$, $\ref{cor-5}$ and $\ref{cor-6}$, we have 
\begin{equation}
S=k(n-1)g.  \label{eq-H-3}
\end{equation}%
By $(\ref{eq-H-2})$ and $(\ref{eq-H-3})$, we get 
\[
E_{1}g(Y,Z)-E_{2}\eta (Y)\eta (Z)=4k(n-1)g(Y,Z) 
\]%
Contracting above equation, we get 
\[
E_{1}n-E_{2}=4kn(n-1). 
\]%
By using the value of $E_{1}$ and $E_{2}$, we get 
\[
c=\frac{4kn^{2}-12nk-3n^{2}+8k+12n-9}{n^{2}-1}. 
\]%
If $M$ is a ${\cal W}_{2}$-semisymmetric $N(k)$-contact metric manifold of
constant $\varphi $-sectional curvature $c$, then by using Corollary $\ref%
{cor-7}$, we have 
\[
c=\frac{4r-8nk-3n^{2}+8k+12n-9}{n^{2}-1}. 
\]
\end{ex}

\section{$({\cal T}_{\!a},S_{{\cal T}_{b}})$-semisymmetry\label{sect-TS}}

\begin{defn-new}
A semi-Riemannian manifold is said to be $({\cal T}_{\!a},S_{{\cal T}_{b}})$%
-semisymmetric if 
\begin{equation}
{\cal T}_{\!a}(V,U)\cdot S_{{\cal T}_{b}}=0.  \label{eq-T.S}
\end{equation}%
In particular, it is said to be $({\cal T}_{\!a},S)$-semisymmetric or, in
short, {\em ${\cal T}_{\!a}$-Ricci-semisymmetric} if it satisfies 
\begin{equation}
{\cal T}_{\!a}(V,U)\cdot S=0,  \label{eq-T.Sa}
\end{equation}%
where ${\cal T}_{\!a}\left( V,U\right) $ acts as a derivation on $S$. In
particular, if in (\ref{eq-T.Sa}), {\em ${\cal T}_{\!a}$} is equal to $R$, $%
{\cal C}_{\ast }$, ${\cal C}$, ${\cal L}$, ${\cal V}$, ${\cal P}_{\ast }$, $%
{\cal P}$, ${\cal M}$, ${\cal W}_{0}$, ${\cal W}_{0}^{\ast }$, ${\cal W}_{1}$%
, ${\cal W}_{1}^{\ast }$, ${\cal W}_{2}$, ${\cal W}_{3}$, ${\cal W}_{4}$, $%
{\cal W}_{5}$, ${\cal W}_{6}$, ${\cal W}_{7}$, ${\cal W}_{8}$, ${\cal W}_{9}$%
, then it becomes Ricci-semisymmetric \cite{Adati-Miyazawa-79}, ${\cal C}%
_{\ast }$-Ricci-semisymmetric, ${\cal C}$-Ricci-semisymmetric (or, Weyl
Ricci-semisymmetric \cite{Ozgur-2005}), ${\cal L}$-Ricci-semisymmetric, $%
{\cal V}$-Ricci-semisymmetric (concircular Ricci-semisymmetric \cite%
{Hong-Ozgur-Tripathi-06}), ${\cal P}_{\ast }$-Ricci-semisymmetric, ${\cal P}$%
-Ricci-semisymmetric, ${\cal M}$-Ricci-semisymmetric, ${\cal W}_{0}$%
-Ricci-semisymmetric, ${\cal W}_{0}^{\ast }$-Ricci-semisymmetric, ${\cal W}%
_{1}$-Ricci-semisymmetric, ${\cal W}_{1}^{\ast }$-Ricci-semisymmetric, $%
{\cal W}_{2}$-Ricci-semisymmetric, ${\cal W}_{3}$-Ricci-semisymmetric, $%
{\cal W}_{4}$-Ricci-semisymmetric, ${\cal W}_{5}$-Ricci-semisymmetric, $%
{\cal W}_{6}$-Ricci-semisymmetric, ${\cal W}_{7}$-Ricci-semisymmetric, $%
{\cal W}_{8}$-Ricci-semisymmetric, ${\cal W}_{9}$-Ricci-semisymmetric,
respectively. 
\end{defn-new}

\begin{lem}
\label{th-T-S} Let $M$ be an $n$-dimensional $({\cal T}_{\!a},S_{{\cal T}%
_{b}})$-semisymmetric $\left( N(k),\xi \right) $-semi-Riemannian manifold.
Then 
\begin{eqnarray}
0 &=&\ \varepsilon a_{5}(b_{0}+nb_{1}+b_{2}+b_{3}+b_{5}+b_{6})S^{2}(Y,U) 
\nonumber \\
&&+\ \left\{ \varepsilon (b_{0}+nb_{1}+b_{2}+b_{3}+b_{5}+b_{6})\times \right.
\nonumber \\
&&(-ka_{0}+k(n-1)a_{1}+k(n-1)a_{2}-a_{7}r)  \nonumber \\
&&\qquad \left. +\ \varepsilon (a_{1}+a_{5})(b_{4}r+(n-1)b_{7}r)\right\}
S(Y,U)  \nonumber \\
&&+\ \left\{ \varepsilon k(n-1)(a_{2}+a_{4})(b_{4}r+(n-1)b_{7}r)\right. 
\nonumber \\
&&\qquad +\ \varepsilon k(n-1)(b_{0}+nb_{1}+b_{2}+b_{3}+b_{5}+b_{6})\times 
\nonumber \\
&&\left. (ka_{0}+k(n-1)a_{4}+a_{7}r)\right\} g(Y,U)  \nonumber \\
&&+\ k(n-1)(a_{1}+a_{2}+2a_{3}+a_{4}+a_{5}+2a_{6})\times  \nonumber \\
&&\left\{ (b_{4}r+(n-1)b_{7}r)\right.  \nonumber \\
&&\left. +\ k(n-1)(b_{0}+nb_{1}+b_{2}+b_{3}+b_{5}+b_{6})\right\} \eta
(Y)\eta (U).  \label{eq-ricci-TS}
\end{eqnarray}%
In particular, if $M$ is $({\cal T}_{\!a},S_{{\cal T}_{a}})$-semisymmetric $%
\left( N(k),\xi \right) $-semi-Riemannian manifold, then 
\begin{eqnarray}
0 &=&\ \varepsilon a_{5}(a_{0}+na_{1}+a_{2}+a_{3}+a_{5}+a_{6})S^{2}(Y,U) 
\nonumber \\
&&+\ \left\{ \varepsilon (a_{0}+na_{1}+a_{2}+a_{3}+a_{5}+a_{6})\times \right.
\nonumber \\
&&(-ka_{0}+k(n-1)a_{1}+k(n-1)a_{2}-a_{7}r)  \nonumber \\
&&\qquad \left. +\ \varepsilon (a_{1}+a_{5})(a_{4}r+(n-1)a_{7}r)\right\}
S(Y,U)  \nonumber \\
&&+\ \left\{ \varepsilon k(n-1)(a_{2}+a_{4})(a_{4}r+(n-1)a_{7}r)\right. 
\nonumber \\
&&\qquad +\ \varepsilon k(n-1)(a_{0}+na_{1}+a_{2}+a_{3}+a_{5}+a_{6})\times 
\nonumber \\
&&\left. (ka_{0}+k(n-1)a_{4}+a_{7}r)\right\} g(Y,U)  \nonumber \\
&&+\ k(n-1)(a_{1}+a_{2}+2a_{3}+a_{4}+a_{5}+2a_{6})\times  \nonumber \\
&&\left\{ (a_{4}r+(n-1)a_{7}r)\right.  \nonumber \\
&&\left. +\ k(n-1)(a_{0}+na_{1}+a_{2}+a_{3}+a_{5}+a_{6})\right\} \eta
(Y)\eta (U).  \label{eq-ricci-TS1}
\end{eqnarray}
\end{lem}

\noindent {\bf Proof.} Let $M$ be an $n$-dimensional $({\cal T}_{\!a},S_{%
{\cal T}_{b}})$-semisymmetric $\left( N(k),\xi \right) $-semi-Riemannian
manifold. Then 
\begin{equation}
\left( {\cal T}_{\!a}(X,Y)\cdot S_{{\cal T}_{\!b}}\right) (U,V)=0.
\label{eq-T-S-1}
\end{equation}%
Taking $X=\xi =V$ in (\ref{eq-T-S-1}), we get 
\[
({\cal T}_{\!a}(\xi ,Y)\cdot S_{{\cal T}_{\!b}})(U,\xi )=0, 
\]%
which gives 
\begin{equation}
S_{{\cal T}_{\!b}}({\cal T}_{\!a}(\xi ,Y)U,\xi )+S_{{\cal T}_{\!b}}(U,{\cal T%
}_{\!a}(\xi ,Y)\xi )=0.  \label{eq-T-S}
\end{equation}%
Using (\ref{eq-cond}), (\ref{eq-ricci}), (\ref{eq-xi-X-xi}), (\ref{eq-xi-Y-Z}%
), (\ref{eq-ric-T1}) and (\ref{eq-ric-T2}) in (\ref{eq-T-S}), we get (\ref%
{eq-ricci-TS}). $\blacksquare $ \medskip

By Lemma~\ref{th-T-S}, we have the following

\begin{th}
\label{cor-RS} Let $M$ be a $(R,S_{{\cal T}_{a}})$-semisymmetric $\left(
N(k),\xi \right) $-semi-Riemannian manifold such that 
\[
a_{0}+na_{1}+a_{2}+a_{3}+a_{5}+a_{6}\neq 0. 
\]%
Then 
\[
S(X,Y)=k(n-1)g(X,Y). 
\]
\end{th}

\begin{th}
An Einstein manifold is $\left( R,S_{{\cal T}_{a}}\right) $-semisymmetric.
\end{th}

\noindent {\bf Proof.} Let $M$ be an Einstein manifold. Then we have 
\[
S(X,Y)=d_{1}g(X,Y), 
\]%
where $d_{1}$ is smooth function on the manifold. By (\ref{eq-ric-T}), we
have 
\begin{eqnarray}
S_{{\cal T}_{\!a}}(X,Y)
&=&d_{1}(a_{0}+na_{1}+a_{2}+a_{3}+na_{4}+a_{5}+a_{6}+n(n-1)a_{7})g(X,Y) 
\nonumber \\
&=&d_{2}g(X,Y),  \label{eq-S-T-a-1}
\end{eqnarray}%
where 
\[
d_{2}=d_{1}(a_{0}+na_{1}+a_{2}+a_{3}+na_{4}+a_{5}+a_{6}+n(n-1)a_{7}). 
\]%
Then by condition $(R(X,Y)\cdot S_{{\cal T}_{a}})(U,V)$ and (\ref{eq-S-T-a-1}%
), we get 
\[
-S_{{\cal T}_{a}}(R(X,Y)U,V)-S_{{\cal T}_{a}}(U,R(X,Y)V)=0. 
\]%
Hence this proves the result. $\blacksquare $

\begin{th}
Let $M$ be an Einstein manifold such that 
\[
{\cal T}_{\!a}\in \left\{ R,{\cal C}_{\ast },{\cal C},{\cal L},{\cal V},%
{\cal M},{\cal W}_{0},{\cal W}_{0}^{\ast },{\cal W}_{3}\right\} . 
\]%
Then it is ${\cal T}_{\!a}$-Ricci-semisymmetric.
\end{th}

\noindent {\bf Proof.} Let $M$ be an Einstein manifold such that 
\[
{\cal T}_{\!a}\in \left\{ R,{\cal C}_{\ast },{\cal C},{\cal L},{\cal V},%
{\cal M},{\cal W}_{0},{\cal W}_{0}^{\ast },{\cal W}_{3}\right\} . 
\]%
Then we have 
\begin{equation}
S=\alpha g,  \label{eq-S-alpha}
\end{equation}%
where $\alpha $ is smooth function on the manifold $M$. Using condition $(%
{\cal T}_{\!a}(X,Y)\cdot S)(U,V)$ and (\ref{eq-S-alpha}), we get 
\[
-S({\cal T}_{\!a}(X,Y)U,V)-S(U,{\cal T}_{\!a}(X,Y)V)=0. 
\]%
Hence this proves the result. $\blacksquare $

By Lemma~\ref{th-T-S}, we have the following

\begin{th}
\label{GCT-rss} Let $M$ be a ${\cal T}_{\!a}$-Ricci-semisymmetric $\left(
N(k),\xi \right) $-semi-Riemannian manifold. Then 
\[
\varepsilon a_{5}\,S^{2}(X,Y)=E\,S(X,Y)+Fg(X,Y)+G\eta (X)\eta \left(
Y\right) , 
\]%
where 
\[
E=(\varepsilon ka_{0}+\varepsilon a_{7}r-\varepsilon k(n-1)a_{1}-\varepsilon
k(n-1)a_{2}), 
\]%
\[
F=-\varepsilon k(n-1)(ka_{0}+k(n-1)a_{4}+a_{7}r), 
\]%
\[
G=-k^{2}(n-1)^{2}(a_{1}+a_{2}+2a_{3}+a_{4}+a_{5}+2a_{6}). 
\]%
In particular, if $a_{5}=0\neq E$, then $M$ is $\eta $-Einstein manifold.
\end{th}

In view of Theorem~\ref{GCT-rss}, we have the following

\begin{cor}
Let $M$ be an $n$-dimensional Ricci-semisymmetric $\left( N(k),\xi \right) $%
-semi-Riemannian manifold. Then we have the following table\/{\rm :}~%
\[
\begin{tabular}{|l|l|}
\hline
${\boldmath M}$ & ${\boldmath S=}$ \\ \hline
$N(k)$-contact metric {\rm \cite{Perrone-92}} & $k(n-1)g$ \\ \hline
Sasakian {\rm \cite{Perrone-92}} & $(n-1)g$ \\ \hline
Kenmotsu {\rm \cite{Hong-Ozgur-Tripathi-06}} & $-\,(n-1)g$ \\ \hline
$(\varepsilon )$-Sasakian & $\varepsilon (n-1)g$ \\ \hline
para-Sasakian {\rm (\cite{Adati-Miyazawa-79}, \cite{Pandey-Verma-99}, \cite%
{Sharfuddin-Deshmukh-Husain-80})} & $-\,(n-1)g$ \\ \hline
$(\varepsilon )$-para-Sasakian {\rm \cite{TKYK-09}} & $-\,\varepsilon (n-1)g$
\\ \hline
\end{tabular}%
\]
\end{cor}

\begin{cor}
Let $M$ be an $n$-dimensional ${\cal C}_{\ast }$-Ricci-semisymmetric $\left(
N(k),\xi \right) $-semi-Riemannian manifold. Then we have the following
table\/{\rm :}~%
\[
\begin{tabular}{|l|l|}
\hline
${\boldmath M}$ & ${\boldmath S^{2}=}$ \\ \hline
$N(k)$-contact metric & $-\,\left( \left( k-\dfrac{r}{n(n-1)}\right) \dfrac{%
a_{0}}{a_{1}}-\dfrac{2r}{n}\right) S$ \\ 
& $+k(n-1)\left( \left( k-\dfrac{r}{n(n-1)}\right) \dfrac{a_{0}}{a_{1}}%
+k(n-1)-\dfrac{2r}{n}\right) g$ \\ \hline
Sasakian & $-\,\left( \left( 1-\dfrac{r}{n(n-1)}\right) \dfrac{a_{0}}{a_{1}}-%
\dfrac{2r}{n}\right) S$ \\ 
& $+(n-1)\left( \left( 1-\dfrac{r}{n(n-1)}\right) \dfrac{a_{0}}{a_{1}}+(n-1)-%
\dfrac{2r}{n}\right) g$ \\ \hline
Kenmotsu & $\left( \left( 1+\dfrac{r}{n(n-1)}\right) \dfrac{a_{0}}{a_{1}}+%
\dfrac{2r}{n}\right) S$ \\ 
& $+(n-1)\left( \left( 1+\dfrac{r}{n(n-1)}\right) \dfrac{a_{0}}{a_{1}}+(n-1)+%
\dfrac{2r}{n}\right) g$ \\ \hline
$(\varepsilon )$-Sasakian & $-\,\varepsilon \left( \left( 1-\dfrac{%
\varepsilon r}{n(n-1)}\right) \dfrac{a_{0}}{a_{1}}-\dfrac{2\varepsilon r}{n}%
\right) S$ \\ 
& $+\varepsilon (n-1)\left( \left( \varepsilon -\dfrac{r}{n(n-1)}\right) 
\dfrac{a_{0}}{a_{1}}+\varepsilon (n-1)-\dfrac{2r}{n}\right) g$ \\ \hline
para-Sasakian & $\left( \left( 1+\dfrac{r}{n(n-1)}\right) \dfrac{a_{0}}{a_{1}%
}+\dfrac{2r}{n}\right) S$ \\ 
& $+(n-1)\left( \left( 1+\dfrac{r}{n(n-1)}\right) \dfrac{a_{0}}{a_{1}}+(n-1)+%
\dfrac{2r}{n}\right) g$ \\ \hline
$(\varepsilon )$-para-Sasakian & $\varepsilon \left( \left( 1+\dfrac{%
\varepsilon r}{n(n-1)}\right) \dfrac{a_{0}}{a_{1}}+\dfrac{2\varepsilon r}{n}%
\right) S$ \\ 
& $+\varepsilon (n-1)\left( \left( \varepsilon +\dfrac{r}{n(n-1)}\right) 
\dfrac{a_{0}}{a_{1}}+\varepsilon (n-1)+\dfrac{2r}{n}\right) g$ \\ \hline
\end{tabular}%
\]
\end{cor}

\begin{cor}
Let $M$ be an $n$-dimensional ${\cal C}$-Ricci-semisymmetric $\left(
N(k),\xi \right) $-semi-Riemannian manifold. Then we have the following
table\/{\rm :}~%
\[
\begin{tabular}{|l|l|}
\hline
${\boldmath M}$ & ${\boldmath S^{2}=}$ \\ \hline
$N(k)$-contact metric & $\left( \dfrac{r}{n-1}+k(n-2)\right) S-k(r-(n-1)k)g$
\\ \hline
Sasakian & $\left( \dfrac{r}{n-1}+(n-2)\right) S-(r-(n-1))g$ \\ \hline
Kenmotsu & $\left( \dfrac{r}{n-1}-(n-2)\right) S+(r+(n-1))g$ \\ \hline
$(\varepsilon )$-Sasakian & $\left( \dfrac{r}{n-1}+\varepsilon (n-2)\right)
S-\varepsilon (r-(n-1)\varepsilon )g$ \\ \hline
para-Sasakian {\rm \cite{Ozgur-2005}} & $\left( \dfrac{r}{n-1}-(n-2)\right)
S+(r+(n-1))g$ \\ \hline
$(\varepsilon )$-para-Sasakian & $\left( \dfrac{r}{n-1}-\varepsilon
(n-2)\right) S+\varepsilon (r+(n-1)\varepsilon )g$ \\ \hline
\end{tabular}%
\]
\end{cor}

\begin{cor}
Let $M$ be an $n$-dimensional ${\cal L}$-Ricci-semisymmetric $\left(
N(k),\xi \right) $-semi-Riemannian manifold. Then we have the following
table\/{\rm :}~%
\[
\begin{tabular}{|l|l|}
\hline
${\boldmath M}$ & ${\boldmath S^{2}=}$ \\ \hline
$N(k)$-contact metric & $k(n-2)S+k^{2}(n-1)g$ \\ \hline
Sasakian & $(n-2)S+(n-1)g$ \\ \hline
Kenmotsu & $-\,(n-2)S+(n-1)g$ \\ \hline
$(\varepsilon )$-Sasakian & $\varepsilon (n-2)S+(n-1)g$ \\ \hline
para-Sasakian & $-\,(n-2)S+(n-1)g$ \\ \hline
$(\varepsilon )$-para-Sasakian & $-\,\varepsilon (n-2)S+(n-1)g$ \\ \hline
\end{tabular}%
\]
\end{cor}

\begin{cor}
Let $M$ be an $n$-dimensional ${\cal V}$-Ricci-semisymmetric $\left(
N(k),\xi \right) $-semi-Riemannian manifold. Then we have the following
table\/{\rm :}~%
\[
\begin{tabular}{|l|l|}
\hline
${\boldmath M}$ & {\rm Result} \\ \hline
$N(k)$-contact metric & $S=k(n-1)g\ {\rm or}\ r=kn(n-1)$ \\ \hline
Sasakian & $S=(n-1)g\ {\rm or}\ r=n(n-1)$ \\ \hline
Kenmotsu {\rm \cite{Hong-Ozgur-Tripathi-06}} & $S=-(n-1)g\ {\rm or}\
r=-n(n-1)$ \\ \hline
$(\varepsilon )$-Sasakian & $S=\varepsilon (n-1)g\ {\rm or}\ r=\varepsilon
n(n-1)$ \\ \hline
para-Sasakian {\rm \cite{Ozgur-Tripathi-2007}} & $S=-(n-1)g\ {\rm or}\
r=-n(n-1)$ \\ \hline
$(\varepsilon )$-para-Sasakian & $S=-\varepsilon (n-1)g\ {\rm or}\
r=-\varepsilon n(n-1)$ \\ \hline
\end{tabular}%
\]
\end{cor}

\begin{cor}
Let $M$ be an $n$-dimensional ${\cal P}_{\ast }$-Ricci-semisymmetric $\left(
N(k),\xi \right) $-semi-Riemannian manifold such that $a_{0}+(n-1)a_{1}\neq
0 $. Then we have the following table\/{\rm :}~%
\[
\begin{tabular}{|l|l|}
\hline
${\boldmath M}$ & {\rm Result} \\ \hline
$N(k)$-contact metric & $S=k(n-1)g\ {\rm or}\ r=\dfrac{n(n-1)ka_{0}}{%
a_{0}+(n-1)a_{1}}$ \\ \hline
Sasakian & $S=(n-1)g\ {\rm or}\ r=\dfrac{n(n-1)a_{0}}{a_{0}+(n-1)a_{1}}$ \\ 
\hline
Kenmotsu & $S=-(n-1)g\ {\rm or}\ r=\dfrac{-n(n-1)a_{0}}{a_{0}+(n-1)a_{1}}$
\\ \hline
$(\varepsilon )$-Sasakian & $S=\varepsilon (n-1)g\ {\rm or}\ r=\dfrac{%
n(n-1)\varepsilon a_{0}}{a_{0}+(n-1)a_{1}}$ \\ \hline
para-Sasakian & $S=-(n-1)g\ {\rm or}\ r=\dfrac{-n(n-1)a_{0}}{a_{0}+(n-1)a_{1}%
}$ \\ \hline
$(\varepsilon )$-para-Sasakian & $S=-\varepsilon (n-1)g\ {\rm or}\ r=\dfrac{%
-n(n-1)\varepsilon a_{0}}{a_{0}+(n-1)a_{1}}$ \\ \hline
\end{tabular}%
\]
\end{cor}

\begin{cor}
Let $M$ be an $n$-dimensional ${\cal P}$-Ricci-semisymmetric $\left(
N(k),\xi \right) $-semi-Riemannian manifold. Then we have the following
table\/{\rm :}~%
\[
\begin{tabular}{|l|l|}
\hline
${\boldmath M}$ & ${\boldmath S=}$ \\ \hline
$N(k)$-contact metric & $k(n-1)g$ \\ \hline
Sasakian & $(n-1)g$ \\ \hline
Kenmotsu & $-\,(n-1)g$ \\ \hline
$(\varepsilon )$-Sasakian & $\varepsilon (n-1)g$ \\ \hline
para-Sasakian & $-\,(n-1)g$ \\ \hline
$(\varepsilon )$-para-Sasakian & $-\,\varepsilon (n-1)g$ \\ \hline
Lorentzian para-Sasakian & $(n-1)g$ \\ \hline
\end{tabular}%
\]
\end{cor}

\begin{cor}
Let $M$ be an $n$-dimensional ${\cal M}$-Ricci-semisymmetric $\left(
N(k),\xi \right) $-semi-Riemannian manifold. Then we have the following
table\/{\rm :}~%
\[
\begin{tabular}{|l|l|}
\hline
${\boldmath M}$ & ${\boldmath S^{2}=}$ \\ \hline
$N(k)$-contact metric & $2k(n-1)S+k^{2}(n-1)^{2}g$ \\ \hline
Sasakian & $2(n-1)S+(n-1)^{2}g$ \\ \hline
Kenmotsu & $-\,2(n-1)S+(n-1)^{2}g$ \\ \hline
$(\varepsilon )$-Sasakian & $2\varepsilon (n-1)S+(n-1)^{2}g$ \\ \hline
para-Sasakian & $-\,2(n-1)S+(n-1)^{2}g$ \\ \hline
$(\varepsilon )$-para-Sasakian & $-2\varepsilon (n-1)S+(n-1)^{2}g$ \\ \hline
\end{tabular}%
\]
\end{cor}

\begin{cor}
Let $M$ be an $n$-dimensional ${\cal W}_{0}$-Ricci-semisymmetric $\left(
N(k),\xi \right) $-semi-Riemannian manifold. Then we have the following
table\/{\rm :}~%
\[
\begin{tabular}{|l|l|}
\hline
${\boldmath M}$ & ${\boldmath S^{2}=}$ \\ \hline
$N(k)$-contact metric & $2k(n-1)S-k^{2}(n-1)^{2}g$ \\ \hline
Sasakian & $2(n-1)S-(n-1)^{2}g$ \\ \hline
Kenmotsu & $-\,2(n-1)S-(n-1)^{2}g$ \\ \hline
$(\varepsilon )$-Sasakian & $2\varepsilon (n-1)S-(n-1)^{2}g$ \\ \hline
para-Sasakian & $-\,2(n-1)S-(n-1)^{2}g$ \\ \hline
$(\varepsilon )$-para-Sasakian & $-\,2\varepsilon (n-1)S-(n-1)^{2}g$ \\ 
\hline
\end{tabular}%
\]
\end{cor}

\begin{cor}
Let $M$ be an $n$-dimensional ${\cal W}_{0}^{\ast }$-Ricci-semisymmetric $%
\left( N(k),\xi \right) $-semi-Riemannian manifold. Then we have the
following table\/{\rm :}~%
\[
\begin{tabular}{|l|l|}
\hline
${\boldmath M}$ & ${\boldmath S^{2}=}$ \\ \hline
$N(k)$-contact metric & $k^{2}(n-1)^{2}g$ \\ \hline
Sasakian & $(n-1)^{2}g$ \\ \hline
Kenmotsu & $(n-1)^{2}g$ \\ \hline
$(\varepsilon )$-Sasakian & $(n-1)^{2}g$ \\ \hline
para-Sasakian & $(n-1)^{2}g$ \\ \hline
$(\varepsilon )$-para-Sasakian & $(n-1)^{2}g$ \\ \hline
\end{tabular}%
\]
\end{cor}

\begin{cor}
Let $M$ be an $n$-dimensional ${\cal W}_{1}$-Ricci-semisymmetric $\left(
N(k),\xi \right) $-semi-Riemannian manifold. Then we have the following
table\/{\rm :}~%
\[
\begin{tabular}{|l|l|}
\hline
${\boldmath M}$ & ${\boldmath S=}$ \\ \hline
$N(k)$-contact metric & $k(n-1)g$ \\ \hline
Sasakian & $(n-1)g$ \\ \hline
Kenmotsu & $-\,(n-1)g$ \\ \hline
$(\varepsilon )$-Sasakian & $\varepsilon (n-1)g$ \\ \hline
para-Sasakian & $-\,(n-1)g$ \\ \hline
$(\varepsilon )$-para-Sasakian & $-\,\varepsilon (n-1)g$ \\ \hline
\end{tabular}%
\]
\end{cor}

\begin{cor}
Let $M$ be an $n$-dimensional ${\cal W}_{1}^{\ast }$-Ricci-semisymmetric $%
\left( N(k),\xi \right) $-semi-Riemannian manifold. Then we have the
following table\/{\rm :}~%
\[
\begin{tabular}{|l|l|}
\hline
${\boldmath M}$ & ${\boldmath S=}$ \\ \hline
$N(k)$-contact metric & $k(n-1)g$ \\ \hline
Sasakian & $(n-1)g$ \\ \hline
Kenmotsu & $-\,(n-1)g$ \\ \hline
$(\varepsilon )$-Sasakian & $\varepsilon (n-1)g$ \\ \hline
para-Sasakian & $-\,(n-1)g$ \\ \hline
$(\varepsilon )$-para-Sasakian & $-\,\varepsilon (n-1)g$ \\ \hline
\end{tabular}%
\]
\end{cor}

\begin{cor}
Let $M$ be an $n$-dimensional ${\cal W}_{2}$-Ricci-semisymmetric $\left(
N(k),\xi \right) $-semi-Riemannian manifold. Then we have the following
table\/{\rm :}~%
\[
\begin{tabular}{|l|l|}
\hline
${\boldmath M}$ & ${\boldmath S^{2}=}$ \\ \hline
$N(k)$-contact metric & $k(n-1)S$ \\ \hline
Sasakian & $(n-1)S$ \\ \hline
Kenmotsu & $-\,(n-1)S$ \\ \hline
$(\varepsilon )$-Sasakian & $\varepsilon (n-1)S$ \\ \hline
para-Sasakian & $-\,(n-1)S$ \\ \hline
$(\varepsilon )$-para-Sasakian & $-\,\varepsilon (n-1)S$ \\ \hline
\end{tabular}%
\]
\end{cor}

\begin{cor}
Let $M$ be an $n$-dimensional ${\cal W}_{3}$-Ricci-semisymmetric $\left(
N(k),\xi \right) $-semi-Riemannian manifold. Then we have the following
table\/{\rm :}~%
\[
\begin{tabular}{|l|l|}
\hline
${\boldmath M}$ & ${\boldmath S=}$ \\ \hline
$N(k)$-contact metric & $k(n-1)g$ \\ \hline
Sasakian & $(n-1)g$ \\ \hline
Kenmotsu & $-\,(n-1)g$ \\ \hline
$(\varepsilon )$-Sasakian & $\varepsilon (n-1)g$ \\ \hline
para-Sasakian & $-\,(n-1)g$ \\ \hline
$(\varepsilon )$-para-Sasakian & $-\,\varepsilon (n-1)g$ \\ \hline
\end{tabular}%
\]
\end{cor}

\begin{cor}
Let $M$ be an $n$-dimensional ${\cal W}_{4}$-Ricci-semisymmetric $\left(
N(k),\xi \right) $-semi-Riemannian manifold. Then we have the following
table\/{\rm :}~%
\[
\begin{tabular}{|l|l|}
\hline
${\boldmath M}$ & ${\boldmath S^{2}=}$ \\ \hline
$N(k)$-contact metric & $k(n-1)S-k^{2}(n-1)^{2}g+\varepsilon
k^{2}(n-1)^{2}\eta \otimes \eta $ \\ \hline
Sasakian & $(n-1)S-(n-1)^{2}g+(n-1)^{2}\eta \otimes \eta $ \\ \hline
Kenmotsu & $-\,(n-1)S-(n-1)^{2}g+(n-1)^{2}\eta \otimes \eta $ \\ \hline
$(\varepsilon )$-Sasakian & $\varepsilon (n-1)S-(n-1)^{2}g+\varepsilon
(n-1)^{2}\eta \otimes \eta $ \\ \hline
para-Sasakian & $-\,(n-1)S-(n-1)^{2}g+(n-1)^{2}\eta \otimes \eta $ \\ \hline
$(\varepsilon )$-para-Sasakian & $-\,\varepsilon
(n-1)S-(n-1)^{2}g+\varepsilon (n-1)^{2}\eta \otimes \eta $ \\ \hline
\end{tabular}%
\]
\end{cor}

\begin{cor}
Let $M$ be an $n$-dimensional ${\cal W}_{5}$-Ricci-semisymmetric $\left(
N(k),\xi \right) $-semi-Riemannian manifold. Then we have the following
table\/{\rm :}~%
\[
\begin{tabular}{|l|l|}
\hline
${\boldmath M}$ & ${\boldmath S^{2}=}$ \\ \hline
$N(k)$-contact metric & $2k(n-1)S-k^{2}(n-1)^{2}g$ \\ \hline
Sasakian & $2(n-1)S-(n-1)^{2}g$ \\ \hline
Kenmotsu & $-\,2(n-1)S-(n-1)^{2}g$ \\ \hline
$(\varepsilon )$-Sasakian & $2\varepsilon (n-1)S-(n-1)^{2}g$ \\ \hline
para-Sasakian & $-\,2(n-1)S-(n-1)^{2}g$ \\ \hline
$(\varepsilon )$-para-Sasakian & $-\,2\varepsilon (n-1)S-(n-1)^{2}g$ \\ 
\hline
\end{tabular}%
\]
\end{cor}

\begin{cor}
Let $M$ be an $n$-dimensional ${\cal W}_{6}$-Ricci-semisymmetric $\left(
N(k),\xi \right) $-semi-Riemannian manifold. Then we have the following
table\/{\rm :}~%
\[
\begin{tabular}{|l|l|}
\hline
${\boldmath M}$ & ${\boldmath2S=}$ \\ \hline
$N(k)$-contact metric & $k(n-1)g+k(n-1)\eta \otimes \eta $ \\ \hline
Sasakian & $(n-1)g+(n-1)\eta \otimes \eta $ \\ \hline
Kenmotsu & $-\,(n-1)g-(n-1)\eta \otimes \eta $ \\ \hline
$(\varepsilon )$-Sasakian & $\varepsilon (n-1)g+(n-1)\eta \otimes \eta $ \\ 
\hline
para-Sasakian & $-\,(n-1)g-(n-1)\eta \otimes \eta $ \\ \hline
$(\varepsilon )$-para-Sasakian & $-\,\varepsilon (n-1)g-(n-1)\eta \otimes
\eta $ \\ \hline
\end{tabular}%
\]
\end{cor}

\begin{cor}
Let $M$ be an $n$-dimensional ${\cal W}_{7}$-Ricci-semisymmetric $\left(
N(k),\xi \right) $-semi-Riemannian manifold. Then we have the following
table\/{\rm :}~%
\[
\begin{tabular}{|l|l|}
\hline
${\boldmath M}$ & ${\boldmath S=}$ \\ \hline
$N(k)$-contact metric & $k(n-1)g$ \\ \hline
Sasakian & $(n-1)g$ \\ \hline
Kenmotsu & $-\,(n-1)g$ \\ \hline
$(\varepsilon )$-Sasakian & $\varepsilon (n-1)g$ \\ \hline
para-Sasakian & $-\,(n-1)g$ \\ \hline
$(\varepsilon )$-para-Sasakian & $-\,\varepsilon (n-1)g$ \\ \hline
\end{tabular}%
\]
\end{cor}

\begin{cor}
Let $M$ be an $n$-dimensional ${\cal W}_{8}$-Ricci-semisymmetric $\left(
N(k),\xi \right) $-semi-Riemannian manifold. Then we have the following
table\/{\rm :}~%
\[
\begin{tabular}{|l|l|}
\hline
${\boldmath M}$ & ${\boldmath2S=}$ \\ \hline
$N(k)$-contact metric & $k(n-1)g+k(n-1)\eta \otimes \eta $ \\ \hline
Sasakian & $(n-1)g+(n-1)\eta \otimes \eta $ \\ \hline
Kenmotsu & $-\,(n-1)g-(n-1)\eta \otimes \eta $ \\ \hline
$(\varepsilon )$-Sasakian & $\varepsilon (n-1)g+(n-1)\eta \otimes \eta $ \\ 
\hline
para-Sasakian & $-\,(n-1)g-(n-1)\eta \otimes \eta $ \\ \hline
$(\varepsilon )$-para-Sasakian & $-\,\varepsilon (n-1)g-(n-1)\eta \otimes
\eta $ \\ \hline
\end{tabular}%
\]
\end{cor}

\begin{cor}
Let $M$ be an $n$-dimensional ${\cal W}_{9}$-Ricci-semisymmetric $\left(
N(k),\xi \right) $-semi-Riemannian manifold. Then we have the following
table\/{\rm :}~%
\[
\begin{tabular}{|l|l|}
\hline
${\boldmath M}$ & ${\boldmath S=}$ \\ \hline
$N(k)$-contact metric & $k(n-1)\eta \otimes \eta $ \\ \hline
Sasakian & $(n-1)\eta \otimes \eta $ \\ \hline
Kenmotsu & $-\,(n-1)\eta \otimes \eta $ \\ \hline
$(\varepsilon )$-Sasakian & $(n-1)\eta \otimes \eta $ \\ \hline
para-Sasakian & $-\,(n-1)\eta \otimes \eta $ \\ \hline
$(\varepsilon )$-para-Sasakian & $-\,(n-1)\eta \otimes \eta $ \\ \hline
\end{tabular}%
\]
\end{cor}

\medskip

\noindent Department of Mathematics\newline
Faculty of Science\newline
Banaras Hindu University\newline
Varanasi-221005\newline
mmtripathi66@yahoo.com \medskip

\noindent Department of Mathematics\newline
Faculty of Science\newline
Banaras Hindu University\newline
Varanasi-221005\newline
punam\_2101@yahoo.co.in

\end{document}